\documentclass[11pt,a4paper,reqno]{amsart}
\usepackage[T1]{fontenc}
\usepackage{lmodern,amsmath,amssymb,amsthm,mathtools,bm,microtype}
\usepackage[margin=27mm]{geometry}
\usepackage{needspace,etoolbox,enumitem,xcolor}
\usepackage[colorlinks=true,linkcolor=blue!50!black,citecolor=blue!50!black,urlcolor=blue!50!black]{hyperref}
\pretocmd{\section}{\Needspace{10\baselineskip}}{}{}
\pretocmd{\subsection}{\Needspace{5\baselineskip}}{}{}
\newtheorem{theorem}{Theorem}[section]
\newtheorem{lemma}[theorem]{Lemma}
\newtheorem{proposition}[theorem]{Proposition}
\newtheorem{corollary}[theorem]{Corollary}
\theoremstyle{remark}
\newtheorem{remark}[theorem]{Remark}
\newcommand{\R}{\mathbb R}

\newcommand{\Q}{\mathbb Q}
\newcommand{\Z}{\mathbb Z}
\newcommand{\T}{\mathbb T}

\newcommand{\dd}{\,\mathrm d}

\newcommand{\LL}{\mathcal L}
\DeclareMathOperator{\dist}{dist}

\numberwithin{equation}{section}
\setlist[enumerate]{itemsep=4pt,topsep=6pt}

\title[Diophantine conditions for absolute continuity of projections]{Diophantine conditions for absolute continuity of projections for self-similar measures}

\author{Xiong Jin}
\address{Department of Mathematics, The University of Manchester,
Oxford Road, Manchester M13 9PL, United Kingdom}
\email{xiong.jin@manchester.ac.uk}
\author{Tuomas Sahlsten}
\address{Department of Mathematics and Statistics, University of Helsinki,
P.O. Box 68, FI-00014 University of Helsinki, Finland}
\email{tuomas.sahlsten@helsinki.fi}



\hypersetup{pdftitle={Diophantine conditions for absolute continuity of projections for self-similar measures},
pdfauthor={Xiong Jin and Tuomas Sahlsten}}

\thanks{T.S. acknowledges support from the Research Council of Finland's Academy Research Fellowship \emph{``Quantum chaos of large and many body systems''}, grant Nos. 347365, 353738. Part of this research was carried out during X.J. and T.S.'s visits to the programme at Institut Mittag-Leffler \emph{Interactions between Fractal Geometry, Harmonic Analysis and Dynamical Systems} (Autumn 2026).}

\begin{document}
\begin{abstract} We prove that every projection of a planar self-similar measure $\mu$ is equivalent to Lebesgue measure on its support if the rotations satisfy a simultaneous Diophantine condition with exponent $\tau$ and $\mu$ has finite $s$-energy for some $1+2\tau\le s<2$. Under strong separation, we prove the fibre measures have exact dimension $\dim_H\mu-1$. For self-similar sets $E$ with $\dim_HE>1+2\tau$, the same rotation condition implies uniformly positive projection lengths and fibres of dimension $\dim_HE-1$ at Lebesgue-almost every point of each image without separation assumptions. The absolute continuity and fibre dimension theorems extend to non-coordinate projections of products under a joint Diophantine condition with exponent $\tau$ on the logarithmic contraction ratios. These imply for all self-similar sets $K,L\subset\mathbb R$ satisfying this condition and $\dim_HK+\dim_HL>1+2\tau$, every sum $uK+L$, $u\ne0$, has positive length, and $\dim_H(L\cap(y-uK))=\dim_HK+\dim_HL-1$ for Leb. a.e. $y\in uK+L$. The proof combines the $L^1$ Gaussian smoothing approach used in the higher dimensional case with averaging at random stopping scales with a Brownian martingale bound using It\^{o} calculus.
\end{abstract}

\maketitle

\section{Introduction}\label{sec:intro}

\subsection{Background} 

Marstrand-Mattila's projection theorem \cite{Mar54} says that for any Borel $E \subset \R^d$ and for almost every plane $V \in G(d,k)$, orthogonal projection $\pi_V E$ onto $V$ preserves Hausdorff dimension when $\dim_H E \leq k$ and has positive volume when $\dim_H E > k$.  A central question is when these extend to \textit{restricted families} of $V$ through curvature or other structure in the family, or to \textit{every} $V$ under additional assumptions on the set or measure, such as self-similarity or non-linearity. Projection and slicing problems have become a highly active area connecting multiple areas from geometric measure theory and harmonic analysis with additive combinatorics, effective equidistribution, to quantum chaos, see e.g. surveys \cite{FFJ15,Fa26}.

In this topic, we recently proved in \cite{JS26} a uniform projection and slicing theorems for self-similar sets in $\R^d$, $d\ge3$, whose iterated function systems (IFSs) have rotations generating a dense subgroup of $\mathrm{SO}(d)$. If $E$ is such a set, then all its $k$-dimensional orthogonal projections have uniformly positive Lebesgue measure whenever $1\le k<\dim_HE$, and for every such projection the fibres have dimension $\dim_HE-k$ at Lebesgue-almost every point of the image. These results require no separation condition and follow from an absolute continuity theorem for projections of self-similar measures. The main idea in \cite{JS26} is to estimate $L^1$ differences between Gaussian smoothings of the projected measures at successive scales. Self-similarity allows these quantities to be averaged over rotations without introducing growing factors under rescaling. This makes it possible to prove absolute continuity using Varj\'u's weak gap theorem \cite[Corollary 7]{Va13} under dense generation alone.

Extending this argument to $\R^2$ is difficult because the rotation group $\mathrm{SO}(2)\cong\T=\R/\Z$ is abelian, whereas Varj\'u's theorem applies to compact connected semisimple Lie groups. On $\T$, dense generation alone is not enough, the rate of mixing depends on the arithmetic properties of the rotation angles. Indeed, Rapaport's counterexample \cite[Theorem 1.1]{Ra17} shows that a planar self-similar measure can have dense rotations and dimension greater than one while some of its projections are singular. If one goes beyond self-similarity, Algom, Rodriguez-Hertz and Wang \cite{ARW26b} proved absolute continuity for nonlinear self-conformal and strongly irreducible, proximal self-affine measures using Littlewood-Paley and renewal methods.

Our aim is to replace Varj\'u's weak gap theorem by estimates using simultaneous Diophantine conditions on the rotations or contractions. These conditions control the concentration of the directions produced by the rotation walk at random stopping scales. Combining these estimates with Gaussian smoothing and a Brownian martingale estimate, we have uniform logarithmic moment bounds for the projected densities under suitable energy hypotheses. In particular, every projection is absolutely continuous. The main novelty of our approach is the use of Brownian motion and It\^o calculus to prove logarithmic moment bounds for projections of \textit{deterministic} self-similar measures by transforming the Gaussian smoothing estimates into martingale estimates, which may be of independent interest.

The same method applies to products of self-similar sets and measures on the line under a corresponding Diophantine condition on the logarithmic contraction ratios. Here the relative contraction scales determine the projection directions, and we prove results for every non-coordinate projection. This implies applications to both sums and intersections: under the hypotheses below, every sum $uK+L$, $u\ne0$, has positive Lebesgue measure, and $L\cap(y-uK)$ has dimension $\dim_HK+\dim_HL-1$ for Lebesgue-almost every $y\in uK+L$. The set results allow arbitrary overlaps in the defining IFSs.

Let us now describe our results in detail.

\subsection{Planar self-similar sets with Diophantine rotations}\label{sec:set-results}

Let $\LL$ denote Lebesgue measure on $\R$. For $\theta\in\T=\R/\Z$, define $\pi_\theta:\R^2\to\R$ by $\pi_\theta x=x\cdot(\cos 2\pi\theta,\sin 2\pi\theta)$. For planar self-similar sets with dense rotations, Peres and Shmerkin \cite{PS09} proved that every orthogonal projection has dimension $\min\{\dim_HE,1\}$. Shmerkin and Solomyak \cite[Theorem C]{SS16} proved that, for any planar self-similar set of dimension greater than one, the directions with zero projection length form a set of Hausdorff dimension zero. Our first result is an explicit rotation condition under which \textit{every} projection has positive length and the expected fibre dimension almost everywhere. Given $\tau>0$, we say $\bm\alpha=(\alpha_1,\ldots,\alpha_q)\in\R^q$, $q\ge1$, satisfies $(D)_\tau$ if there is $c>0$ such that
\begin{equation}\label{eq:dio-int}
 \max_{1\le j\le q}\|k\alpha_j\|_{\T}\ge c k^{-\tau}
 \qquad\text{for all }k=1,2,\ldots, \qquad \|t\|_{\T}=\dist(t,\Z).
\end{equation}

\begin{theorem}\label{thm:pl-sets}
Let $E\subset\R^2$ be the attractor of an IFS $\mathcal S=(S_i)_{i=1}^m$, where
$$
 S_i(x)=r_iO_{\alpha_i}x+t_i,\qquad
 0<r_i<1,\quad t_i\in\R^2,
$$
and $O_\alpha$ is rotation by $2\pi\alpha$. Suppose that, for some $\tau>0$,
\begin{equation}\label{eq:pl-dio}
 \bm\alpha(\mathcal S):=(\alpha_1,\ldots,\alpha_m)
 \quad\text{satisfies }(D)_\tau.
\end{equation}
If $\dim_HE>1+2\tau$, then
$$
 \inf_{\theta\in\T}\LL(\pi_\theta E)>0.
$$
Moreover, for every fixed $\theta\in\T$,
$$
 \dim_H\bigl(E\cap\pi_\theta^{-1}\{y\}\bigr)=\dim_HE-1
 \quad\text{for }\LL\text{-almost every }y\in\pi_\theta E.
$$
\end{theorem}

\begin{remark}\label{rem:tau-range}
\begin{itemize}
    \item[(1)] By Dirichlet's simultaneous approximation theorem, we have $\tau\ge1/m$. The badly approximable case corresponds to $\tau=1/m$. The hypotheses of Theorem \ref{thm:pl-sets} require $\frac1m\le\tau<\frac{\dim_HE-1}{2}\le\frac12,$
so $m(\dim_HE-1)>2$ and $m\ge3$. The larger $m$ is, the smaller $\tau$ can be, and therefore the closer the dimension threshold $\dim_HE>1+2\tau$ can be to the sharp bound $\dim_HE>1$. For $\tau>1/m$, the simultaneous Jarn\'ik-Besicovitch theorem (see \cite{BRV16}) implies that the set of vectors $\bm\alpha\in\R^m$ which do not satisfy $(D)_\tau$ has Hausdorff dimension $\frac{m+1}{1+\tau}<m$, see Proposition \ref{prop:dio-exc}. In particular, $(D)_\tau$ holds for Lebesgue-almost every $\bm\alpha\in\R^m$.
\item[(2)] Under dense rotations and the open set condition, Falconer and the first author \cite[Theorem 4.8]{FJ15} proved weak dimension conservation: outside a zero-dimensional set of directions, fibres of dimension at least $\dim_HE-1-\varepsilon$ occur over a set of positive length for every $\varepsilon\in(0,\dim_HE-1)$. Rapaport \cite{Ra20} proved uniformly positive projection lengths and the fibre dimension formula in every direction for homogeneous systems $z\mapsto\lambda z+t_i$ satisfying the strong separation condition (SSC), with attractor dimension greater than one and $\arg\lambda/\pi\notin\Q$, provided $\lambda$ avoids a fixed zero-dimensional exceptional set. Theorem \ref{thm:pl-sets} implies an explicit sufficient condition without exceptional sets and allows different contraction ratios and rotations. Explicit examples with algebraic rotation angles are given in Subsection \ref{sec:pl-ex}.
\item[(3)] We note that Diophantine conditions on rotations have also been used in quantitative projection estimates. Ero\u{g}lu \cite[Theorem 3.1]{Ero07} proved quantitative Favard-length decay for a class of homogeneous self-similar systems under a Diophantine rotation condition and additional structural assumptions. Also, recently Shmerkin and Usuki \cite[Theorem 1.5]{SU26} proved that every affine slice of a planar self-similar set $E$ has upper box dimension at most $\max\{\dim_HE-1,0\}$ under the open set condition and a Diophantine condition.
\end{itemize}
\end{remark}

\subsection{Products, sums and intersections of self-similar sets} 

Our second result concerns products $K\times L$ and imposes a Diophantine condition on the contraction ratios. Let $K,L\subset\R$ be the attractors of IFSs $\mathcal S=(S_i)_{i=1}^m$ and $\mathcal S'=(S'_j)_{j=1}^n$, respectively, where $S_i(x)=r_ix+d_i$, $S'_j(x)=r'_jx+d'_j$, $0<|r_i|,|r'_j|<1.$ Peres and Shmerkin \cite[Theorem 2]{PS09} proved that $\dim_H(uK+L)=\min\{\dim_HK+\dim_HL,1\}$ for every $u\ne0$ if $\log|r_i|/\log|r_j'|\notin\Q$ for some $i,j$. When $\dim_HK+\dim_HL>1$, Shmerkin \cite[Theorem 1.3]{Sh14} proved positive length of $uK+L$ for every $u\ne0$ for homogeneous systems with SSC and common contraction ratios $r$ for $K$ and $s$ for $L$, provided $\log r/\log s\notin\Q$ and $r$ avoids a fixed zero-dimensional exceptional set. Shmerkin and Solomyak \cite[Theorem E]{SS16} proved that, for arbitrary self-similar sets $K,L\subset\R$ with $\dim_HK+\dim_HL>1$, the set of scalings $u\ne0$ for which $\LL(uK+L)=0$ has Hausdorff dimension zero. 

For intersections, Shmerkin \cite[Corollary 7.3]{Sh19a} proved $\overline{\dim}_B(L\cap(y-uK))\le\max\{\dim_HK+\dim_HL-1,0\}$ for all $u,y$ when the sets are homogeneous and strongly separated under irrationality condition. Shmerkin's work and Wu \cite{Wu19} independently resolved Furstenberg's intersection conjecture. Recently, Shmerkin and Usuki \cite[Theorem 1.1]{SU26} extended slicing bounds to non-homogeneous systems under the open set condition and a Diophantine condition. See also a different direction to these problems using thickness and relative position \cite{HKY93}.

Now under a Diophantine condition on the logarithmic contraction ratios from both IFSs, we can avoid non-explicit exceptional sets and prove a projection theorem that holds for \textit{all} projections except the coordinate directions. Given $\tau>0$, we define that $\bm\alpha=(\alpha_1,\ldots,\alpha_q)\in\R^q$ satisfies $(D)'_\tau$ if there is $c>0$ such that
\begin{equation}\label{eq:dio}
 \max_{1\le j\le q}\|t\alpha_j\|_{\T}\ge c t^{-\tau}
 \qquad\text{for all }t\ge1.
\end{equation}
Note that Lemma \ref{lem:real-freq}, $(\alpha_1,\ldots,\alpha_q)$ satisfies $(D)_\tau$ if and only if $(1,\alpha_1,\ldots,\alpha_q)$ satisfies $(D)'_\tau$. Our result under this condition is:

\begin{theorem}\label{thm:sums}
 Let $K,L\subset\R$ be the attractors of IFSs $\mathcal S=(S_i)_{i=1}^m$ and $\mathcal S'=(S'_j)_{j=1}^n$ for $S_i(x)=r_ix+d_i$, $S'_j(x)=r'_jx+d'_j$ and $0<|r_i|,|r'_j|<1$. Write $\beta_i=\log(1/|r_i|)$ and $\beta'_j=\log(1/|r'_j|)$. Suppose that, for some $\tau>0$,
\begin{equation}\label{eq:sum-dio}
 \bm\beta(\mathcal S,\mathcal S'):=
 (\beta_1,\ldots,\beta_m,\beta'_1,\ldots,\beta'_n)
 \quad\text{satisfies }(D)'_\tau.
\end{equation}
If $\dim_HK+\dim_HL>1+2\tau$, then $uK+L$ has positive Lebesgue measure for every $u\ne0$, with
$$
 \inf_{u\in J}\LL(uK+L)>0
$$
for every nonempty compact $J\subset\R\setminus\{0\}$. Moreover, for every fixed $u\ne0$,
$$
 \dim_H\bigl(L\cap(y-uK)\bigr)=\dim_HK+\dim_HL-1
 \quad\text{for }\LL\text{-almost every }y\in uK+L.
$$
\end{theorem}

\begin{remark}\label{rem:sum-tau}
The vector in \eqref{eq:sum-dio} has $m+n$ positive entries. By Lemma \ref{lem:real-freq}, dividing by one entry reduces \eqref{eq:sum-dio} to condition \eqref{eq:dio-int} for the remaining $m+n-1$ ratios. By Dirichlet's simultaneous approximation theorem, we have $\tau\ge1/(m+n-1)$. If the vector
\begin{equation}\label{eq:one-fac}
 \bm\beta(\mathcal S)
 :=(\beta_1,\ldots,\beta_m)
\end{equation}
satisfies $(D)'_\tau$, then \eqref{eq:sum-dio} holds for every IFS $\mathcal{S}'$, since adjoining entries cannot decrease the maximum in \eqref{eq:dio}. The joint condition also allows examples in which neither IFS separately has a suitable Diophantine exponent, see Subsection \ref{sec:log-ex}.

For two homogeneous systems, the combined logarithmic vector has at most two distinct entries. By Lemma \ref{lem:real-freq} and Dirichlet's theorem, $(D)'_\tau$ then requires $\tau\ge1$, which contradicts $2\ge\dim_HK+\dim_HL>1+2\tau$. Thus Theorem \ref{thm:sums} does \textit{not} apply to a pair of homogeneous defining IFSs.
\end{remark}

The constructions in Subsection \ref{sec:log-ex} give the following consequence, in which one fixed set $K$ works for every self-similar set $L$ above the stated dimension threshold.

\begin{corollary}\label{cor:fix}
For every $\varepsilon\in(0,1)$, there is an explicitly defined self-similar set $K\subset\R$ with $\dim_HK=1-\varepsilon/2$ such that, for every self-similar set $L\subset\R$ with $\dim_HL>\varepsilon$, the sumset $uK+L$ has positive Lebesgue measure for every $u\ne0$, with
$$
 \inf_{u\in J}\LL(uK+L)>0
$$
for every nonempty compact $J\subset\R\setminus\{0\}$. Moreover, for each fixed $u\ne0$,
$$
 \dim_H\bigl(L\cap(y-uK)\bigr)=\dim_HL-\varepsilon/2
$$
for $\LL$-almost every $y\in uK+L$.
\end{corollary}

The product $K\times L$ is a diagonal self-affine set. Related projection results for self-affine carpets and strongly irreducible planar self-affine sets appear e.g. in \cite{FJS10,BKK26}, and slicing estimates for carpets and other planar self-affine sets are given in \cite{Alg20,AW23,BKY26}.

For the proofs, we use a common formulation of the planar and product results. We write $(\mathrm I)_\tau$ for the planar systems satisfying \eqref{eq:pl-dio}, and $(\mathrm{II})_\tau$ for pairs of line systems satisfying \eqref{eq:sum-dio}. Let $E$ denote the planar attractor in $(\mathrm I)$ and the product $K\times L$ in $(\mathrm{II})$. The admissible projection vectors is the set:
$$
\mathcal D=
\begin{cases}
 \R^2\setminus\{0\},&\text{in }(\mathrm I),\\
 \{v\in\R^2:v_1v_2\ne0\},&\text{in }(\mathrm{II}).
\end{cases}
$$
For $v\in\R^2$, write $\pi_vx=v\cdot x$ and $E_v=\pi_vE$. Up to a nonzero scalar multiple of the image, we may take $v=(\cos 2\pi\theta,\sin 2\pi\theta)$ in $(\mathrm I)$ and $v=(u,1)$, $u\ne0$, in $(\mathrm{II})$.

\begin{theorem}\label{thm:sets}
Let $E$ be the planar self-similar set in $(\mathrm{I})_\tau$, or the product of self-similar sets $K\times L$ in $(\mathrm{II})_\tau$. If $\dim_HE>1+2\tau$, then for every compact subset $\mathcal V\subset\mathcal D$,
$$
 \inf_{v\in\mathcal V}\LL(E_v)>0,
$$
and, for every fixed $v\in\mathcal D$,
$$
 \dim_H\bigl(E\cap \pi_v^{-1}\{y\}\bigr)=\dim_HE-1
 \quad\text{for }\LL\text{-almost every }y\in E_v.
$$
\end{theorem}

Theorem \ref{thm:pl-sets} follows from Theorem \ref{thm:sets} by taking $v=(\cos 2\pi\theta,\sin 2\pi\theta)$. Theorem \ref{thm:sums} follows from Theorem \ref{thm:sets} by using the identity $\pi_{(u,1)}(K\times L)=uK+L$ and the similarity $y'\mapsto((y-y')/u,y')$ from $L\cap(y-uK)$ onto $(K\times L)\cap \pi_{(u,1)}^{-1}\{y\}$, together with $\dim_H(K\times L)=\dim_HK+\dim_HL$ (see \cite[Corollary 7.4]{Fa03}).

For every $\varepsilon\in(0,1)$, we construct in Subsection \ref{sec:log-ex} an explicit self-similar set $K\subset\R$ satisfying SSC, with $\bm\beta(\mathcal S)$ in \eqref{eq:one-fac} satisfying $(D)'_{\varepsilon/4}$ and with $\dim_HK=1-\varepsilon/2$. Thus the dimension hypothesis in Theorem \ref{thm:sets} and the Diophantine condition \eqref{eq:sum-dio} hold for every self-similar set $L\subset\R$ with $\dim_HL>\varepsilon$, and we prove Corollary \ref{cor:fix}.

\subsection{Results on measures}\label{sec:measure-results}

Our set theoretical results will be consequences on analogous absolute continuity and dimension conservation theorems for self-similar measures, which we will now state separately. For an IFS $(S_i)_{i\in\mathcal I}$ and a positive probability vector $p=(p_i)_{i\in\mathcal I}$, the associated \textit{self-similar measure} is the unique probability measure $\mu$ satisfying
\begin{equation}\label{eq:ss}
 \mu=\sum_{i\in\mathcal I}p_i(S_i)_*\mu.
\end{equation}
The IFS with attractor $E$ has strong separation condition (SSC) if $S_i(E)\cap S_j(E)=\varnothing$, $i\ne j$. 
 
 Fix $\tau>0$. We say that $\mu$ belongs to $(\mathrm I)_\tau$ if it is a planar self-similar measure with positive weights whose defining IFS has rotation angles $2\pi\alpha_1,\ldots,2\pi\alpha_m$ such that $(\alpha_1,\ldots,\alpha_m)$ satisfies $(D)_\tau$ in \eqref{eq:dio-int}. We say that $\mu$ belongs to $(\mathrm{II})_\tau$ if $\mu=\nu\otimes\nu'$, where $\nu,\nu'$ are self-similar measures on $\R$ with positive weights and contraction multipliers $r_1,\ldots,r_m$ and $r'_1,\ldots,r'_n$, respectively, such that $(\log(1/|r_1|),\ldots,\log(1/|r_m|),\log(1/|r'_1|),\ldots,\log(1/|r'_n|))$ satisfies $(D)'_\tau$ in \eqref{eq:dio}. These Diophantine conditions concern the defining IFSs, no arithmetic condition is imposed on the weights. For $v\in\mathcal D$, write $\mu_v=(\pi_v)_*\mu$.

For each fixed $v\in\mathcal D$, let $(\mu_{v,y})_{y\in\R}$ be a disintegration of $\mu$ over $\mu_v$, so $\mu=\int\mu_{v,y}\dd\mu_v(y)$. The conditional measure $\mu_{v,y}$ is a probability measure supported on the fibre $E\cap \pi_v^{-1}\{y\}$ for $\mu_v$-almost every $y$. For a probability measure $\lambda$, write $\dim_H\lambda$ for the upper Hausdorff dimension of $\lambda$, that is, the infimum of $\dim_HB$ over Borel sets $B$ of full $\lambda$-measure. We call $\lambda$ exact dimensional if $
 \lim_{r\downarrow0}\frac{\log\lambda(B(x,r))}{\log r}=\dim_H\lambda$ for $\lambda\text{-almost every }x.$
For $0<s<2$, the $s$-energy of $\lambda$ is $
 I_s(\lambda)=\iint|x-z|^{-s}\dd\lambda(x)\dd\lambda(z).$
The measure version of Marstrand's projection theorem implies $L^2$ densities in almost every direction for compactly supported planar probability measures with finite $s$-energy for some $1\le s<2$, see \cite[Theorem 9.7]{Ma95}. Our Diophantine hypothesis yields absolute continuity in every admissible direction, together with a uniform logarithmic moment bound.

The following measure theorem implies Theorem \ref{thm:sets} and hence Theorems \ref{thm:pl-sets} and \ref{thm:sums}.

\begin{theorem}\label{thm:ac}
Let $\mu$ belong to $(\mathrm{I})_\tau$ or $(\mathrm{II})_\tau$, and suppose that $I_s(\mu)<\infty$ for some $1+2\tau\le s<2$. Then, for every $v\in\mathcal D$, there is a density $f_v$ such that
$$
 \mu_v=f_v\dd x\sim\LL|_{E_v}.
$$
Here $\sim$ denotes equivalence of measures. Write $\kappa=(s-1)/(2\tau)\ge1$. For every compact $\mathcal V\subset\mathcal D$, there is $C_{\mathcal V}<\infty$ such that
\begin{equation}\label{eq:log-mom}
 \sup_{v\in\mathcal V}\int_{\R}
 f_v(x)\bigl(\log(e+f_v(x))\bigr)^\kappa\dd x
 \le C_{\mathcal V}.
\end{equation}
In particular, $\inf_{v\in\mathcal V}\LL(E_v)>0$. If SSC holds in $(\mathrm{I})$, or in both IFSs in $(\mathrm{II})$, then, for each fixed $v\in\mathcal D$, the conditional measures $\mu_{v,y}$ are exact dimensional of dimension $\dim_H\mu-1$ for $\LL$-almost every $y\in E_v$.
\end{theorem}

For $u\ne0$, let $\nu_u$ be the pushforward of $\nu$ under $y\mapsto uy$. Since $\mu_{(u,1)}=\nu_u*\nu'$ and $E_{(u,1)}=uK+L$, Theorem \ref{thm:ac} implies the following convolution statement.

\begin{corollary}\label{cor:conv}
Let $\nu$ and $\nu'$ be as in $(\mathrm{II})_\tau$, and suppose $I_s(\nu\otimes\nu')<\infty$ for some $1+2\tau\le s<2$. Then, for every $u\ne0$,$
 \nu_u*\nu' \sim\LL|_{uK+L}.$
\end{corollary}

For projection dimensions, Hochman and Shmerkin \cite[Theorem 1.6]{HS12} proved the maximal-dimension formula under dense rotations and SSC, and Falconer and the first author \cite{FJ14} removed the separation assumption. Shmerkin and Usuki \cite[Theorems 1.1 and 1.5]{SU26} proved $L^q$-dimension formulas for projections and convolutions in non-homogeneous settings under the open set condition and their Diophantine hypotheses.

For absolute continuity, Shmerkin and Solomyak \cite{SS16} proved that the exceptional directions for a homogeneous planar self-similar measure of dimension greater than one satisfying SSC have Hausdorff dimension zero. K\"aenm\"aki and Orponen \cite[Proposition 1.7]{KO23} proved an exceptional-direction result for homotheties with unequal contraction ratios, assuming similarity dimension greater than one, non-collinear fixed points and exponential separation.

Related projection-dimension results for self-conformal Gibbs measures, self-conformal measures and diagonal self-affine measures appear in \cite{BJ19,Py25,Py25b}. Algom, Rodriguez Hertz and Wang \cite[Theorems 1.2 and 1.3]{ARW26b} give criteria for absolute continuity of projections of planar self-conformal measures whose IFS is not analytically conjugate to a self-similar IFS, and for strongly irreducible proximal self-affine measures satisfying suitable energy and Furstenberg-measure dimension hypotheses. This does not apply to our self-similar dense rotation setting no our self-affine product systems as the linear parts are reducible.

To prove the set results from the measure results, we approximate the attractors by SSC subsystems while preserving the Diophantine exponent. This follows the approach as in the higher dimensional case \cite{JS26}, see Section \ref{sec:sets}.

\subsection*{Structure of the paper}

Section \ref{sec:strategy} outlines the proof, and Section \ref{sec:prelim} defines the symbolic notation. Section \ref{sec:avg} proves the Brownian martingale estimate and reduces the logarithmic moment bound to an average over stopped projection vectors. Section \ref{sec:sob} proves the directional interval estimates and completes the logarithmic moment bound using the energy of the measure. Section \ref{sec:ac} proves absolute continuity, equivalence with Lebesgue measure, and dimension conservation, and proves the set results by approximation with strongly separated subsystems. Finally, in the Appendix \ref{sec:ex} we construct explicit examples for planar systems and products satisfying our Diophantine conditions.

\section{Proof strategy}\label{sec:strategy}

The proof of main measure theoretic Theorem \ref{thm:ac} that implies the set theoretical results is based on a uniform logarithmic moment bound for Gaussian smoothings of the projected measures. A Brownian martingale estimate reduces this to a weighted second moment integral. Self-similarity then allows us to average over projection directions at random stopping scales. The Diophantine conditions control the concentration of these directions, and integration in scale implies a bound in terms of the energy $I_s(\mu)$. These estimates replace Varj\'u's weak gap theorem used in the earlier higher dimensional version by the authors \cite{JS26}. 

Our use of lazy walks is motivated by Borda's quantitative convergence estimates on $\T$ \cite[Theorem 6]{Bor21}. The use of Diophantine conditions on contraction logarithms and renewal walks is also related to the Fourier decay argument of Li and the second author \cite{LS22}.

Let us describe the strategy more in details and compare to the higher dimensional case in \cite{JS26} in the end. We use the symbolic notation introduced in Section \ref{sec:ifs}.

\subsection{The logarithmic moment bound}\label{sec:strategy-moments}

Our goal is a bound that is uniform in both the smoothing scale and the projection directiobn. This allows us to pass to the limit as the smoothing scale tends to zero.

Let $\mu$ belong to $(\mathrm I)_\tau$ or $(\mathrm{II})_\tau$, with $I_s(\mu)<\infty$ for some $1+2\tau\le s<2$, recall the formulation of Theorem \ref{thm:ac}. For a direction $v\in\mathcal{D}$ and $h>0$ define the Gaussian smoothing of $\mu_v$ by
$$
 f_{v,h}(x)=(\mu_v*\Gamma_h)(x)
 =\int_\R
 \Gamma_h(x-y)\dd\mu_v(y),\quad x\in\mathbb{R},
$$
where $\Gamma_h(x)=\frac{1}{\sqrt{2\pi}h}e^{-x^2/(2h^2)}$ is the centred Gaussian density with standard deviation $h$. Fix a compact set $\mathcal V\subset\mathcal D$ and write $\kappa=\frac{s-1}{2\tau}\ge1.$ The main estimate we want to prove is:
\begin{equation}\label{eq:ren-ent}
 \sup_{v\in\mathcal V}\sup_{0<h<h_0}
 \int_\R f_{v,h}\bigl(\log^+(h_0f_{v,h})\bigr)^\kappa\dd x
 \le C_0\bigl(1+I_s(\mu)\bigr)
\end{equation}
for fixed constants $C_0<\infty$ and $h_0\in(0,1)$, depending only on the IFSs, their probability vectors, $s$, $\tau$ and $\mathcal V$, see Proposition \ref{prop:bm-ren}. This implies uniform integrability and hence absolute continuity.

\subsection{Reduction to Brownian Martingales and using It\^o calculus}\label{sec:strategy-brownian}

To prove \eqref{eq:ren-ent}, we first convert the problem into an estimate for spatial derivatives of the smoothings. Here the crucial feature is the \textit{Gaussian semigroup property} that allows us to construct a Brownian martingale whose quadratic variation contains these derivatives, where we can apply \textit{It\^o calculus}.

Firstly, we have the advertised semigroup property for the Gaussian kernel:
\begin{equation}\label{eq:g-sem}
 \Gamma_h*\Gamma_a=\Gamma_{\sqrt{h^2+a^2}}\text{ for } h,a>0.
\end{equation}
We use this to construct a Brownian martingale. Let $(B_t)_{t\ge0}$ be standard Brownian motion. Fix $v\in\mathcal V$ and $0<\varepsilon<h_0\le1$. Let $T=h_0^2-\varepsilon^2$ and $h_t=\sqrt{h_0^2-t}$. For $x\in\mathbb{R}$ define
$$
 M_t^x=f_{v,h_t}(x+B_t),\qquad 0\le t\le T.
$$
Let $(\mathcal F_t)_{t\ge0}$ be the Brownian filtration. For $0\le t<T$, by independence of Brownian increments, \eqref{eq:g-sem} and $\varepsilon^2+T-t=h_t^2$, we have:
\begin{equation}\label{eq:bm-cond}
 \mathbb E\bigl[f_{v,\varepsilon}(x+B_T)\mid\mathcal F_t\bigr]=\bigl(f_{v,\varepsilon}*\Gamma_{\sqrt{T-t}}\bigr)(x+B_t)=M_t^x.
\end{equation}
Thus $(M_t^x)_{0\le t\le T}$ is a martingale. Then by applying It\^o's formula to it, we have:
$$
 \dd M_t^x=f_{v,h_t}'(x+B_t)\dd B_t,
$$
and its quadratic variation process is $\dd\langle M^x\rangle_t=|f'_{v,h_t}(x+B_t)|^2\dd t$. By the Gaussian bounds of $f_{v,h}$ in \eqref{eq:bm-bd}, $M^x$ is bounded and $M_0^x\le h_0^{-1}$. For the logarithmic moment estimate, define the running maximum
$$
 \overline M_t^x=h_0^{-1}\vee\sup_{0\le r\le t}M_r^x,
 \qquad 0\le t\le T.
$$
Then by Lemma \ref{lem:bm-ent}, we have:
\begin{equation}\label{eq:bm-level}
 \mathbb E\bigl[M_T^x(\log^+(h_0M_T^x))^\kappa\bigr]
 \le C_\kappa h_0(M_0^x)^2+C_\kappa\mathbb E\int_0^T
 \frac{(1+\log(h_0\overline M_t^x))^{\kappa-1}}{\overline M_t^x}
 \dd\langle M^x\rangle_t.
\end{equation}
Such estimates are often referred to as \textit{good-$\lambda$ inequalities}, introduced by Burkholder and Gundy \cite{BG70}, see also Burkholder \cite{Bur73}. Since $h_0\overline M_t^x\le h_0/h_t$ and $\kappa\ge1$, by \eqref{eq:g-ent}, to prove \eqref{eq:ren-ent}, it is enough to show that, for a fixed $C_1>0$,
\begin{equation}\label{eq:bm-ren}
 \mathcal{A}_{v,\varepsilon}:=\int_\R\mathbb E\int_0^T
 \left(1+\log\frac{h_0}{h_t}\right)^{\kappa-1}
 \frac{|f'_{v,h_t}(x+B_t)|^2}{\overline M_t^x}\dd t\dd x
 \le C_1I_s(\mu).
\end{equation}
We note a related approach using It\^o calculus and Burkholder-Davis-Gundy inequality were used by Fraser and the second author \cite{FS18} to find the Fourier dimension of the graph of $B_t$.

\subsection{Using self-similarity at stopping scales}\label{sec:strategy-stopping}

We next use self-similarity to compare the fine smoothing with a coarser smoothing in a new projection direction. For $H > 0$ define
\begin{equation}\label{eq:h-scale}
 h(H):=e^{-H}H^{-1/(2\tau)}
\end{equation}
The reason for this choice is that we stop at the scale $e^{-H}$ and then we go slightly finer scale to $h(H)$ where we study the smoothing of the projected measure. Let $H_0,H_\varepsilon$ satisfy $h(H_0)=h_0$ and $h(H_\varepsilon)=\varepsilon$. We choose $h_0$ sufficiently small that $H_0\ge\max\{e,1/(2\tau)\}$. Write $t(H)$ for the Brownian time satisfying $h_{t(H)}=h(H)$. Set
\begin{equation}\label{eq:dv-def}
 D_v(H)=e^{-2H}\int_\R\mathbb E
 \frac{|f'_{v,h(H)}(x+B_{t(H)})|^2}
      {\overline M_{t(H)}^x}\dd x.
\end{equation}
Changing variables from $t$ to $H$ and using $\kappa\ge1$, we have (with a constant $C_\kappa'>0$ depending only on $\kappa$) that:
\begin{equation}\label{eq:bm-scale}
 \mathcal{A}_{v,\varepsilon}
 \le C_\kappa'\int_{H_0}^{H_\varepsilon}H^{\kappa-1-1/\tau}D_v(H)\dd H.
\end{equation}

We next bound $D_v(H)$ using a stopped \textit{residual vector}. Let the Bernoulli sequence $\omega$ have law $\mathbb P_p$ and, in $(\mathrm{II})$, let $\omega'$ have law $\mathbb P_{p'}$, independently of $\omega$. Write $\bm\omega=\omega$ in $(\mathrm I)$ and $\bm\omega=(\omega,\omega')$ in $(\mathrm{II})$. For $\omega=(\omega_j)_{j\ge1}$, write $W_n=\sum_{j=1}^n\beta_{\omega_j}$, with $W_0=0$ and $\beta_i=-\log|r_i|$. For $H>0$, let $N_H=\min\{n\ge1:W_n\ge H\}$ and $\omega_H=\omega|_{N_H}$. In $(\mathrm{II})$, let $\omega'_H$ denote the second word stopped when its contraction ratio first becomes at most $e^{-W_{N_H}}$. Define the \textit{residual vector}:
\begin{equation}\label{eq:st-par}
 v_H(\bm\omega)=
 \begin{cases}
  e^H r_{\omega_H}O_{-\alpha_{\omega_H}}v,&(\mathrm I),\\
  e^H(r_{\omega_H}v_1,r'_{\omega'_H}v_2),&(\mathrm{II}).
 \end{cases}
\end{equation}
This is the projection vector after rescaling the stopped cylinder by $e^H$. Since the contraction ratios are bounded away from zero, $v_H(\bm\omega)$ lies in a fixed compact set $\widetilde{\mathcal V}\subset\mathcal D$ for every $v\in\mathcal V$, see \eqref{eq:ang-comp} and \eqref{eq:prod-comp}.

Let $t_{v,H}(\bm\omega)$ be the projection of the translation of the stopped map, and denote by $\Phi_{b,\rho}(x)=b+\rho x$ the affine map. Conditional on the stopped words, the remaining tails have their original Bernoulli laws. Hence
\begin{equation}\label{eq:st-law}
 \mu_v=\mathbb E\left[
 (\Phi_{t_{v,H}(\bm\omega),e^{-H}})_*\mu_{v_H(\bm\omega)}\right].
\end{equation}
Define the following second moment observable
\begin{equation}\label{eq:res-der}
 F_H(v')=\|f'_{v',H^{-1/(2\tau)}}\|_2^2,\quad v'\in\mathcal{D}.
\end{equation}
Here and below, both $L^1$ and $L^2$ norms are with respect to Lebesgue measure. We use \eqref{eq:g-recip} with $a=e^{-H}$, then apply Gaussian scaling and Cauchy-Schwarz to \eqref{eq:st-law}. By \eqref{eq:weighted-bd}, we prove with a fixed $C_2>0$ that:
\begin{equation}\label{eq:bm-stop}
 D_v(H)\le C_2\mathbb E[F_H(v_H(\bm\omega))].
\end{equation}

\subsection{Random stopping scales and the direction law}\label{sec:strategy-directions}

The remaining issue is how the stopped projection directions are distributed. Randomising the stopping scale  implies a law with the interval estimate \eqref{eq:dir-mass}, Lemma \ref{lem:bm-avg} allows us to use this law in the scale integral.

We randomise the stopping scale using \textit{lazy renewal walks} following Borda's approach \cite[Theorem 6]{Bor21}. We prove the interval bound at scale $H^{-1/(2\tau)}$ directly for the angular walk in $(\mathrm I)$. In $(\mathrm{II})$, we prove local bounds at fixed steps and sum over the second walk's steps.

Let $\beta_+=\max_i\beta_i$ and $b_\ast=\max\{1,2\beta_+\}$. Independently of the coding sequences, let $(\mathsf B_j)_{j\ge1}$ be independent Bernoulli variables with parameter $1/2$, and define
$$
 \mathsf K_n=\sum_{j=1}^n\mathsf B_j,\qquad
 n_H=\lfloor H/b_\ast\rfloor,\qquad \mathsf N_H=\mathsf K_{n_H}.
$$
The lazy renewal walk $\widetilde W_n=W_{\mathsf K_n}$ holds with probability $1/2$ and proceeds with probability $1/2$. Define
$$
 Y_H=W_{\mathsf N_H+1} \text{ and } \widetilde v_H=v_{Y_H}.
$$
For $\theta\in\T$, let $e_\theta=(\cos 2\pi\theta,\sin 2\pi\theta)$. We write $\widetilde v_H$ in polar coordinates as
\begin{equation}\label{eq:quot-par}
 \widetilde v_H=\rho_H e_{\vartheta_H}.
\end{equation}
Let $\sigma_{v,H}$ denote the law of $\vartheta_H\in\T$. We prove in Proposition \ref{prop:st-en} that
\begin{equation}\label{eq:dir-mass}
 \sigma_{v,H}(J)\le C_3H^{-1/(2\tau)}
\end{equation}
for $H\ge H_0$, $v\in\mathcal V$ and every arc $J\subset\T$ of length $H^{-1/(2\tau)}$. Thus, at the scale $H^{-1/(2\tau)}$, the direction law of the randomised residual vector is uniformly bounded above by the Haar measure. This is the Diophantine substitute for the rotation averaging estimate in \cite{JS26}.

By \eqref{eq:bm-stop} and Lemma \ref{lem:bm-avg} we obtain
\begin{equation}\label{eq:bm-avg}
 \mathcal{A}_{v,\varepsilon}\le
 C_4\int_{H_0}^\infty H^{\kappa-1-1/\tau}
       \mathbb E[F_H(\widetilde v_H)]\dd H.
\end{equation}

\subsection{The energy estimate}\label{sec:strategy-energy}

The interval estimate controls the directional average in \eqref{eq:bm-avg}. We now integrate over $H$, the choice $\kappa=(s-1)/(2\tau)$ produces the Fourier weight $|\xi|^{s-2}$ associated with the $s$-energy.

By compactness, $\rho_H$ is uniformly bounded above and away from zero, so radial rescaling changes only the constants in the estimates for $F_H$. Using the interval bound \eqref{eq:dir-mass}, an angular variation estimate for $\widehat\mu(R e_\theta)$ and Plancherel's identity, we obtain
$$
\mathbb E[F_H(\widetilde v_H)]\le C_5\int_\R R^2e^{-c_5H^{-1/\tau}R^2}
   \int_\T\bigl(|\widehat\mu(R e_\theta)|^2+H^{-1/\tau}R^2|\nabla\widehat\mu(R e_\theta)|^2\bigr)
       \dd\theta\dd R.
$$
Integrating in $H$ and using polar coordinates, we obtain in Lemma \ref{lem:g-trace} that
\begin{equation}\label{eq:g-trace-int}
 \int_{H_0}^\infty H^{\kappa-1-1/\tau}
       \mathbb E[F_H(\widetilde v_H)]\dd H\le C_6\int_{\R^2}|\xi|^{s-2}
       \bigl(|\widehat\mu(\xi)|^2+|\nabla\widehat\mu(\xi)|^2\bigr)\dd\xi\le C_6'I_s(\mu).
\end{equation}
Together with \eqref{eq:bm-avg}, this proves \eqref{eq:bm-ren} and hence \eqref{eq:ren-ent}.

This completes the logarithmic moment estimate, uniformly in the smoothing scale and in $v\in\mathcal V$. Section \ref{sec:ac} uses this bound to prove equivalence with Lebesgue measure and dimension conservation, and then implies the set results by approximation with strongly separated subsystems.

\begin{remark}\label{rem:dir-mix}
Let us compare with the higher dimensional approach that we introduced in the earlier paper \cite[Section~4.2]{JS26}. Denote
$\chi=-\log\rho>0 $, $R=Be^{-n\chi} \ge 1$ and
$\alpha=(s-1)/2$. Let
$U_B(V)=\|\mu_V*(\Gamma_{1/B}^V-\Gamma_{2/B}^V)\|_1$
be the $L^1$ Gaussian smoothing difference in \cite{JS26},
 and write
$$
 \mathcal Q_R(\theta)
 =\|\mu_{e_\theta}*(\Gamma_{1/R}-\Gamma_{2/R})\|_2^2.
$$
For one-dimensional projections, the estimate in \cite{JS26} is
$$ U_B(V)\lesssim R^{-\alpha}+Re^{-cn^{1/3}}$$
and its counterpart in the planar case is
$$
 \int_{\mathbb T}\mathcal Q_R\,d\sigma_{v,n\chi} \lesssim R^{-2\alpha}
   +R^{1-2\alpha}(n\chi)^{-1/(2\tau)}.
$$
By using the Brownian martingale, It\^o's formula and the good-$\lambda$ estimate
for quadratic variation with a weight depending on the running maximum $\overline{M}_t^x$, we can work directly with quadratic
quantities and keep the full energy exponent $2\alpha$. In \cite{JS26}, the passage to $L^1$ norm by Cauchy-Schwarz reduces this exponent to $\alpha$. For the second term, the Fourier representation of the
Gaussian smoothings implies
$$
 \|\mathcal Q_R\|_1
 +R^{-1}\|\partial_\theta\mathcal Q_R\|_1
 \lesssim R^{-2\alpha}.
$$
The non-concentration estimates (Proposition \ref{prop:st-en}) then bound the angular average $\int_{\mathbb T}\mathcal Q_R\,d\sigma_{v,n\chi}$ by
$\|\mathcal Q_R\|_1
 +(n\chi)^{-1/(2\tau)}\|\partial_\theta\mathcal Q_R\|_1$.
Thus the second term has the same energy factor $R^{-2\alpha}$.

In \cite{JS26} we sum the $L^1$ smoothing differences over dyadic scales. Here instead we bound the logarithmic moments by $I_s(\mu)$ via the weighted quadratic-variation integral for Brownian martingales. For $0<\tau<1/2$, a direct use of Borda's Wasserstein mixing
estimate \cite[Theorem~6]{Bor21}, following the approach of
\cite{JS26}, would
lead to the weaker sufficient condition
$s>1+\frac{4\tau}{1-2\tau}$.
On the other hand, optimising a bound of the form
$R^{-a}+R^b e^{-cn^{1/3}}$, with fixed $a,b>0$ and
$n\asymp\log B$, still produces a stretched exponential
with power $1/3$. Thus using this approach in \cite{JS26} would only change the
constants, but not the $1/3$ exponent, which comes from the iteration of Varj\'u's weak-gap estimate in \cite[Proposition~3.3]{JS26}.
\end{remark}

\section{Logarithmic moment bounds using Brownian martingales and residues}\label{sec:avg}

We begin the proof of Theorem \ref{thm:ac} by reducing the logarithmic moment bound to an estimate involving residue projection vectors. In Subsection \ref{sec:ifs} we introduce the symbolic notations we need. In Subsection \ref{sec:bm-ent}, a Brownian martingale argument bounds the logarithmic moments by the weighted integral $\mathcal A_{v,\varepsilon}$. Subsection \ref{sec:stop} expresses self-similarity at the contraction scale corresponding to the Gaussian smoothing. Finally, Subsection \ref{sec:bm-stop} uses this decomposition to bound $\mathcal A_{v,\varepsilon}$ by an integral of $\mathbb E[F_H(v_H)]$.

\subsection{Symbolic notations}\label{sec:prelim}\label{sec:ifs}

We define here the symbolic notation used in Section \ref{sec:strategy} and the proofs. Let $(S_i)_{i\in\mathcal I}$ be one of the IFSs under consideration, with attractor $A$, positive probability vector $p=(p_i)_{i\in\mathcal I}$, and associated self-similar measure $\lambda$. Write $\Sigma=\mathcal I^{\mathbb N}$ and $\Sigma^*=\bigcup_{n\ge0}\mathcal I^n$, where $\mathcal I^0=\{\varnothing\}$. For a finite word $w=w_1\cdots w_n$, let $|w|=n$, and let $w^-$ denote the word obtained by deleting its last letter when $w\ne\varnothing$. For $\omega\in\Sigma$, write $\omega|_n=\omega_1\cdots\omega_n$ and
$
 [w]=\{\omega\in\Sigma:\omega|_{|w|}=w\}.
$
We set $|\varnothing|=0$, $\omega|_0=\varnothing$, and $[\varnothing]=\Sigma$. We use $w\omega$ for concatenation. For a finite word $w=w_1\cdots w_n$, define
$
 S_w=S_{w_1}\circ\cdots\circ S_{w_n}$, $
 r_w=\prod_{j=1}^n r_{w_j}$ and $
 p_w=\prod_{j=1}^n p_{w_j},
$
where $S_\varnothing=\mathrm{Id}$ and $r_\varnothing=p_\varnothing=1$. In the planar case, write $\alpha_w=\sum_{j=1}^n\alpha_{w_j}\pmod1$, with $\alpha_\varnothing=0$, so that $S_w$ has linear part $r_wO_{\alpha_w}$. For maps on the line, $r_w$ may be negative, and the contraction ratio of $S_w$ is $|r_w|$. 

The coding map $\Pi:\Sigma\to A$ is defined by $
 \Pi(\omega)=\lim_{n\to\infty}S_{\omega|_n}(0)$.
The map $\Pi$ is continuous and onto, and satisfies $\Pi(w\omega)=S_w(\Pi(\omega))$ and $\Pi([w])=S_w(A)$. Moreover, the Bernoulli measure $\mathbb P_p$ is defined by:
$
 \mathbb P_p=\left(\sum_{i\in\mathcal I}p_i\delta_i\right)^{\otimes\mathbb N}
$ which satisfies $\lambda=\Pi_*\mathbb P_p$ and $\mathbb P_p([w])=p_w$. Conditional on a prefix $w$, the remaining sequence has law $\mathbb P_p$. A finite family $\mathcal W\subset\Sigma^*$ is a \textit{partition} if its cylinders are pairwise disjoint and cover $\Sigma$. Iterating the self-similarity identity we have $\sum_{w\in\mathcal W}p_w=1$ and
\begin{equation}\label{eq:part}
 \lambda=\sum_{w\in\mathcal W}p_w(S_w)_*\lambda.
\end{equation}
Under SSC, $\Pi$ is one-to-one and
$\lambda|_{S_w(A)}=p_w(S_w)_*\lambda$ and 
 $\lambda(S_w(A))=p_w$.

In $(\mathrm I)_\tau$, we apply this notation to $\lambda=\mu$. In $(\mathrm{II})_\tau$, we use it for $\nu$ and $\nu'$, with weights $p$ and $p'$ and primed notation for the second IFS. The two coding sequences are independent. If $\mathcal P$ is a finite family of word pairs whose rectangles $[w]\times[w']$ partition $\Sigma\times\Sigma'$, then
$
 \sum_{(w,w')\in\mathcal P}p_wp'_{w'}=1$ and $
 \mu=\sum_{(w,w')\in\mathcal P}p_wp'_{w'}(S_w\times S'_{w'})_*\mu.
$
Conditional on a rectangle $[w]\times[w']$, the remaining sequences are independent and have laws $\mathbb P_p$ and $\mathbb P_{p'}$.

All logarithms are natural. We write $\log^+t=\max\{\log t,0\}$ for $t>0$, with $\log^+0=0$.

\subsection{The Brownian logarithmic moment estimate}\label{sec:bm-ent}

We first recall the notation used in this subsection. Let $\mu$ belong to $(\mathrm I)_\tau$ or $(\mathrm{II})_\tau$, let $1+2\tau\le s<2$, and write $\kappa=(s-1)/(2\tau)\ge1$. Fix $v\in\mathcal D$ and $0<\varepsilon<h_0\le1$. We write $\mu_v=(\pi_v)_*\mu$, where $\pi_v(x)=v\cdot x$, and
$$
 f_{v,h}(x)=\int_\R\Gamma_h(x-y)\dd\mu_v(y),
 \qquad
 \Gamma_h(z)=\frac1{\sqrt{2\pi}h}e^{-z^2/(2h^2)}.
$$

Let $(B_t)_{t\ge0}$ be standard Brownian motion starting at zero, and let $\mathcal F_t$ denote the information generated by its path up to time $t$. Set
$$
 T=h_0^2-\varepsilon^2,
 \qquad
 h_t=\sqrt{h_0^2-t},
 \qquad 0\le t\le T.
$$
Thus the smoothing scale decreases from $h_0$ to $\varepsilon$ as $t$ increases from $0$ to $T$. For each $x\in\R$, define
$$
 M_t^x=f_{v,h_t}(x+B_t),
 \qquad
 \overline M_t^x=h_0^{-1}\vee\sup_{0\le r\le t}M_r^x.
$$
We write $M^x$ for the whole process $(M_t^x)_{0\le t\le T}$. Its initial and terminal values are
$$
 M_0^x=f_{v,h_0}(x),
 \qquad
 M_T^x=f_{v,\varepsilon}(x+B_T).
$$
The maximum $\overline M_t^x$ is then the largest value reached up to time $t$, with the base value $h_0^{-1}$ so that $\log(h_0\overline M_t^x)\ge0$.

The Gaussian semigroup identity $\Gamma_h*\Gamma_a=\Gamma_{\sqrt{h^2+a^2}}$ and independence of Brownian increments imply, for $0\le t<T$,
$$
 \mathbb E[M_T^x\mid\mathcal F_t] =\bigl(f_{v,\varepsilon}*\Gamma_{\sqrt{T-t}}\bigr)(x+B_t) =f_{v,h_t}(x+B_t)=M_t^x,
$$
since $\varepsilon^2+T-t=h_t^2$. This is the martingale property: conditional on the path up to time $t$, the expected terminal value equals the current value. The martingale property can also be verified using heat equation: Set $u(t,y)=f_{v,h_t}(y)$. Since $\frac{\dd}{\dd t}h_t^2=-1$, we have $\partial_tu+\frac12\partial_y^2u=0$ since $f_{v,\sqrt q}$ solves the Gaussian heat equation $\frac{\partial}{\partial q}f_{v,\sqrt q}(y)=\frac12\partial_y^2f_{v,\sqrt q}(y).$ Thus the drift term in It\^o's formula \cite[Chapter 5]{LG16} vanishes so it implies:
$$
 \dd M_t^x=f'_{v,h_t}(x+B_t)\dd B_t,
 \qquad
 \langle M^x\rangle_t
 =\int_0^t|f'_{v,h_r}(x+B_r)|^2\dd r.
$$
Here $\langle M^x\rangle$ denotes the quadratic variation, obtained as the limit in probability of sums of squared increments over time partitions whose mesh tends to zero. The formula above expresses it in terms of the spatial derivatives of the Gaussian smoothings.

We will bound the logarithmic moment of $M_T^x$ by stopping the martingale at successive height levels. Integrating the resulting estimates over the levels produces a weight involving the maximum. After integration in $x$, this reduces the logarithmic moment estimate to the quantity
$$
 \mathcal A_{v,\varepsilon}
 :=\int_\R\mathbb E\int_0^T
 \left(1+\log\frac{h_0}{h_t}\right)^{\kappa-1}
 \frac{|f'_{v,h_t}(x+B_t)|^2}{\overline M_t^x}
 \dd t\dd x.
$$

\begin{lemma}\label{lem:bm-ent}
With the notation above, there is a constant $C_\kappa$, depending only on $\kappa$, such that, for every $x\in\R$,
$$
 \mathbb E\bigl[M_T^x(\log^+(h_0M_T^x))^\kappa\bigr]
 \le C_\kappa h_0(M_0^x)^2
 +C_\kappa\mathbb E\int_0^T
 \frac{(1+\log(h_0\overline M_t^x))^{\kappa-1}}
      {\overline M_t^x}
 \dd\langle M^x\rangle_t
$$
and
\begin{equation}\label{eq:g-ent}
 \int_\R f_{v,\varepsilon}
       \bigl(\log^+(h_0f_{v,\varepsilon})\bigr)^\kappa\dd x
 \le C_\kappa+C_\kappa\mathcal A_{v,\varepsilon}.
\end{equation}
\end{lemma}

\begin{proof}
Recall that $M_0^x=f_{v,h_0}(x)$. Since $\mu_v$ is a probability measure,
$$
 f_{v,h}(x)=\int_{\R}\Gamma_h(x-y)\dd\mu_v(y),
 \qquad
 \Gamma_h(z)=\frac1{\sqrt{2\pi}h}e^{-z^2/(2h^2)},
$$
and hence $\|f_{v,h}\|_\infty\le\|\Gamma_h\|_\infty=(\sqrt{2\pi}h)^{-1}$. Moreover, differentiating under the integral and setting $u=z/h$ implies
$$
 \|f'_{v,h}\|_\infty\le\|\Gamma'_h\|_\infty
 =h^{-2}\sup_{u\in\R}\frac{|u|e^{-u^2/2}}{\sqrt{2\pi}}
 \le\frac1{\sqrt{2\pi}h^2}.
$$
Therefore, for every $h\ge\varepsilon$,
\begin{equation}\label{eq:bm-bd}
 M_0^x\le\frac1{\sqrt{2\pi}h_0}<h_0^{-1},\qquad
 \|f_{v,h}\|_\infty\le C_1\varepsilon^{-1},\qquad
 \|f'_{v,h}\|_\infty\le C_1\varepsilon^{-2},
 \qquad C_1=(2\pi)^{-1/2}.
\end{equation}
Thus, by the identity $M_t^x=\mathbb E[M_T^x\mid\mathcal F_t]$, $M^x_t$ is a bounded, continuous, nonnegative martingale. Fix $x\in\R$. For $\lambda>h_0^{-1}$, let
$$
 \sigma_\lambda^x=\inf\{t\in[0,T]:M_t^x\ge\lambda\}
$$
be its first hitting time, with $\sigma_\lambda^x=\infty$ if the level is not reached. By continuity, $M_{\sigma_\lambda^x}^x=\lambda$ on $\{\sigma_\lambda^x\le T\}$ and $\{\sigma_\lambda^x>t\}=\{\overline M_t^x<\lambda\}$.

We first estimate the probability that the martingale reaches a given level $\lambda$. For this, we stop the martingale at its first hitting time and use optional sampling.

In stochastic analysis \textit{optional sampling} extends the martingale property from deterministic times to stopping times, see e.g. \cite[Chapter 3]{LG16}. Since $M^x$ is bounded and continuous, for the stopping time $S=T\wedge\sigma_\lambda^x$ by optional sampling theorem we have
$$
 \mathbb E[M_T^x\mid\mathcal F_S]=M_S^x,
$$
where the stopped $\sigma$-algebra $\mathcal F_S$ represents the information available up to time $S$. The event $A=\{\sigma_\lambda^x\le T\}=\{\overline M_T^x\ge\lambda\}$ belongs to $\mathcal F_S$, and $M_S^x=\lambda$ on $A$. Multiplying the conditional expectation identity by $\mathbf1_A$ and taking expectations implies:
$$
 \mathbb E[M_T^x;\overline M_T^x\ge\lambda]
 =\mathbb E[M_{T\wedge\sigma_\lambda^x}^x;\sigma_\lambda^x\le T]
 =\lambda\mathbb P(\overline M_T^x\ge\lambda),
$$
where $\mathbb E[X;A]$ denotes $\mathbb E[X\mathbf1_A]$.

We now apply It\^o's formula to the process $(M^x)^2$ at $T\wedge\sigma_\lambda^x$, which gives:
\begin{equation}\label{eq:bm-stop-mom}
 \lambda^2\mathbb P(\overline M_T^x\ge\lambda)
 \le\mathbb E(M_{T\wedge\sigma_\lambda^x}^x)^2
 =(M_0^x)^2+\mathbb E\int_0^T
       \mathbf1_{\{\overline M_t^x<\lambda\}}\dd\langle M^x\rangle_t.
\end{equation}

We next integrate this estimate over $\lambda$ to obtain the logarithmic moment of $M_T^x$. The integral over the levels will give a weight involving the running maximum $\overline M_t^x$ in the quadratic variation term.

Since $M_T^x\le\overline M_T^x$, we integrate over the levels $\lambda>h_0^{-1}$ so that
\[
 \begin{split}
 \mathbb E\bigl[M_T^x(\log^+(h_0M_T^x))^\kappa\bigr]
 &\le\mathbb E\bigl[M_T^x(\log(h_0\overline M_T^x))^\kappa\bigr]\\
 &=\kappa\int_{h_0^{-1}}^\infty
       \frac{(\log(h_0\lambda))^{\kappa-1}}{\lambda}
       \mathbb E[M_T^x;\overline M_T^x\ge\lambda]\dd\lambda\\
 &=\kappa\int_{h_0^{-1}}^\infty
       (\log(h_0\lambda))^{\kappa-1}
       \mathbb P(\overline M_T^x\ge\lambda)\dd\lambda.
 \end{split}
\]
For $a\ge h_0^{-1}$, define
$$
 g_{\kappa,h_0}(a)
 =\kappa\int_a^\infty
       \frac{(\log(h_0\lambda))^{\kappa-1}}{\lambda^2}\dd\lambda.
$$
By \eqref{eq:bm-stop-mom} and Fubini, we obtain
\begin{equation}\label{eq:level-ent}
 \mathbb E\bigl[M_T^x(\log^+(h_0M_T^x))^\kappa\bigr]
 \le g_{\kappa,h_0}(h_0^{-1})(M_0^x)^2
 +\mathbb E\int_0^T
       g_{\kappa,h_0}(\overline M_t^x)\dd\langle M^x\rangle_t.
\end{equation}
With $\lambda=ae^z$, we have
$$
 g_{\kappa,h_0}(a)
 =\frac\kappa a\int_0^\infty
       (\log(h_0a)+z)^{\kappa-1}e^{-z}\dd z.
$$
In particular, $g_{\kappa,h_0}(h_0^{-1})=h_0\Gamma(\kappa+1)$. Since $\kappa\ge1$ and $\log(h_0a)+z\le(1+\log(h_0a))(1+z)$, it follows that
\begin{equation}\label{eq:level-kern-bd}
 g_{\kappa,h_0}(a)
 \le\frac{C_\kappa}{a}(1+\log(h_0a))^{\kappa-1}.
\end{equation}
This proves the first claim of Lemma \ref{lem:bm-ent}. Since $h_r\ge h_t$ for $r\le t$,
$$
 h_0\overline M_t^x
 \le\max\left\{1,\frac{h_0}{\sqrt{2\pi}\,h_t}\right\}
 \le\frac{h_0}{h_t}.
$$
By \eqref{eq:level-kern-bd},
$$
 g_{\kappa,h_0}(\overline M_t^x)
 \le\frac{C_\kappa}{\overline M_t^x}
       \left(1+\log\frac{h_0}{h_t}\right)^{\kappa-1}.
$$

Finally, we integrate in $x$ to obtain the required estimate for the smoothed density. Translation invariance removes the Brownian shift from the terminal value, while the initial term is bounded uniformly using the Gaussian convolution.

We integrate \eqref{eq:level-ent} in $x$. By translation invariance,
$$
 \begin{aligned}
 \int_\R M_T^x\bigl(\log^+(h_0M_T^x)\bigr)^\kappa\dd x
 &=\int_\R f_{v,\varepsilon}(x+B_T)
       \bigl(\log^+(h_0f_{v,\varepsilon}(x+B_T))\bigr)^\kappa\dd x\\
 &=\int_\R f_{v,\varepsilon}(y)
       \bigl(\log^+(h_0f_{v,\varepsilon}(y))\bigr)^\kappa\dd y.
 \end{aligned}
$$
We also have
$$
 h_0\int_\R(M_0^x)^2\dd x
 =h_0\|f_{v,h_0}\|_2^2\le\frac1{2\sqrt\pi},
$$
and hence \eqref{eq:g-ent} follows.
\end{proof}

\subsection{The smoothing scale and stopped residual vector}\label{sec:stop}

Let $\mu$ and $s$ be as in Lemma \ref{lem:bm-ent}, and $\kappa=(s-1)/(2\tau)\ge1$. Fix a compact set $\mathcal V\subset\mathcal D$. Choose $H_0\ge\max\{e,1/(2\tau)\}$, and let $h_0=h(H_0)$, where
$$
 h(H)=e^{-H}H^{-1/(2\tau)},\qquad H>0.
$$
We first change from Brownian time to the scale $H$. This will allow us to use the stopping partitions at scale $e^{-H}$.

Fix $v\in\mathcal V$ and $0<\varepsilon<h_0$. Since $h$ decreases strictly to zero, there is a unique $H_\varepsilon>H_0$ with $h(H_\varepsilon)=\varepsilon$. Set
\begin{equation}\label{eq:bm-time}
 t(H)=h_0^2-h(H)^2,\qquad
 \dd t=2h(H)^2\left(1+\frac1{2\tau H}\right)\dd H.
\end{equation}
Then $h_{t(H)}=h(H)$, $t(H_0)=0$ and $t(H_\varepsilon)=T$.
For $H_0\le H\le H_\varepsilon$, define
$$
 D_v(H)=e^{-2H}\int_\R\mathbb E
 \frac{|f'_{v,h(H)}(x+B_{t(H)})|^2}
      {\overline M_{t(H)}^x}\dd x.
$$
By the change of variables \eqref{eq:bm-time}, we have:
\[
 \begin{split}
 \mathcal{A}_{v,\varepsilon}
 &=2\int_{H_0}^{H_\varepsilon}
 \left(1+H-H_0+\frac{\log(H/H_0)}{2\tau}\right)^{\kappa-1}
 \left(1+\frac1{2\tau H}\right)H^{-1/\tau}D_v(H)\dd H\\
 &\le C_\kappa'\int_{H_0}^{H_\varepsilon}H^{\kappa-1-1/\tau}D_v(H)\dd H.
 \end{split}
\]
Here we use $1+\log(h_0/h(H))=1+H-H_0+\log(H/H_0)/(2\tau)\le2H$ and $\kappa\ge1$, and take $C_\kappa'=2^{\kappa+1}$. This proves
\begin{equation}\label{eq:bm-scale-proof}
 \mathcal A_{v,\varepsilon}
 \le C_\kappa'\int_{H_0}^{H_\varepsilon}
 H^{\kappa-1-1/\tau}D_v(H)\dd H.
\end{equation}

We next define the stopping partitions at scale $e^{-H}$. In both cases, we stop the contraction walks when the accumulated logarithmic contractions reach the required scale.

We now define the \textit{residual vector} in terms of the renewal walks. Write $\bm\omega=\omega$ in $(\mathrm I)$ and $\bm\omega=(\omega,\omega')$ in $(\mathrm{II})$.
In $(\mathrm I)$, let $\beta_i=-\log r_i$ and $\beta_+=\max_i\beta_i$. For the coding sequence $\omega$ with law $\mathbb P_p$, write
$$
 W_n=\sum_{j=1}^n\beta_{\omega_j},\qquad
 \Theta_n=\sum_{j=1}^n\alpha_{\omega_j}\pmod1,\qquad
 W_0=\Theta_0=0.
$$
For $H>0$, let
$$
 N_H=\min\{n\ge1:W_n\ge H\},\qquad
 \mathcal P_H=\{w:r_w\le e^{-H}<r_{w^-}\}.
$$
Thus $\omega_H=\omega|_{N_H}$. Let $Z_H=W_{N_H}-H\in[0,\beta_+)$. Then define the residual vector:
$$
 v_H(\bm\omega):=e^Hr_{\omega_H}O_{-\alpha_{\omega_H}}v=e^{-Z_H}O_{-\Theta_{N_H}}v.
$$
Hence $v_H(\bm\omega)$ belongs to the compact set
\begin{equation}\label{eq:ang-comp}
 \widetilde{\mathcal V}
 =\{\rho O_\theta v:v\in\mathcal V,\
       e^{-\beta_+}\le\rho\le1,\ \theta\in\T\}
 \subset\mathcal D.
\end{equation}

In $(\mathrm{II})$, we stop the first walk at scale $e^{-H}$ as before, but stop the second walk when it reaches the contraction scale of the first stopped word. This ensures that the two contractions are comparable.

Let $\beta_i=-\log|r_i|$ and $\beta_+=\max_i\beta_i$, and use the same definitions of $W_n$, $N_H$ and $Z_H$. Let $\beta'_j=-\log|r'_j|$ and $\beta'_+=\max_j\beta'_j$. For the independent coding sequence $\omega'$ with law $\mathbb P_{p'}$, let $W'_0=0$ and $W'_n=\sum_{j=1}^n\beta'_{\omega'_j}$. For $H>0$, let
$$
 N'_H=\min\{n\ge1:W'_n\ge W_{N_H}\}.
$$
We use the stopped words $\omega_H=\omega|_{N_H}$ and $\omega'_H=\omega'|_{N'_H}$, which define the partition
$$
 \mathcal P_H=\{(w,w'):|r_w|\le e^{-H}<|r_{w^-}|,\
                    |r'_{w'}|\le|r_w|<|r'_{(w')^-}|\}.
$$
For each first stopped word, the second threshold is fixed. These product cylinders therefore form a stopping partition with weights $p_wp'_{w'}$ and independent Bernoulli tails.
Let $Z'_H=W'_{N'_H}-W_{N_H}\in[0,\beta'_+)$. With $\varepsilon_H=\operatorname{sgn}r_{\omega_H}$ and $\varepsilon'_H=\operatorname{sgn}r'_{\omega'_H}$, we define the residual vector as:
$$
 v_H(\bm\omega)
 :=e^H(r_{\omega_H}v_1,r'_{\omega'_H}v_2)
 =(\varepsilon_Hv_1e^{-Z_H},\
   \varepsilon'_Hv_2e^{-Z_H-Z'_H}).
$$
Thus $v_H(\bm\omega)$ belongs to the compact set
\begin{equation}\label{eq:prod-comp}
 \widetilde{\mathcal V}
 =\{(\varepsilon\rho v_1,\varepsilon'\rho'v_2):
      v\in\mathcal V,\ \varepsilon,\varepsilon'\in\{-1,1\},
      \ e^{-\beta_+}\le\rho\le1,\
      e^{-\beta_+-\beta'_+}\le\rho'\le1\}
 \subset\mathcal D.
\end{equation}

We now use these stopping partitions to decompose $\mu_v$ into rescaled copies of projected measures. The residual vectors remain in the compact set $\widetilde{\mathcal V}$, so the estimates for these measures will be uniform.

For $\mathbf w\in\mathcal P_H$, write
$$
 (\Psi_{\mathbf w},m_{\mathbf w})=
 \begin{cases}
  (S_w,p_w),&\mathbf w=w\text{ in }(\mathrm I),\\
  (S_w\times S'_{w'},p_wp'_{w'}),&\mathbf w=(w,w')\text{ in }(\mathrm{II}),
 \end{cases}
$$
and let $t_{v,\mathbf w}=\pi_v(\Psi_{\mathbf w}(0))$. Write $\Phi_{b,\rho}(x)=b+\rho x$ for the affine map.
On the symbolic cylinder corresponding to $\mathbf w$, we have
$$
 t_{v,H}(\bm\omega)=t_{v,\mathbf w},\qquad
 v_H(\bm\omega)=v_{H,\mathbf w},
$$
where $v_{H,\mathbf w}$ is the constant residual vector on this cylinder.
Then
\begin{equation}\label{eq:st-map}
 \pi_v\circ\Psi_{\mathbf w}
 =\Phi_{t_{v,\mathbf w},e^{-H}}\circ\pi_{v_{H,\mathbf w}}.
\end{equation}
The symbolic cylinder corresponding to $\mathbf w$ has Bernoulli probability $m_{\mathbf w}$. Thus, for every nonnegative Borel function $F$,
$$
 \mathbb E[F(v_H)]
 =\sum_{\mathbf w\in\mathcal P_H}
 m_{\mathbf w}F(v_{H,\mathbf w}).
$$
Applying $\pi_v$ to \eqref{eq:part}, we have $\sum_{\mathbf w\in\mathcal P_H}m_{\mathbf w}=1$ and
\begin{equation}\label{eq:st-mix}
 \mu_v=\sum_{\mathbf w\in\mathcal P_H}m_{\mathbf w}
 (\Phi_{t_{v,\mathbf w},e^{-H}})_*\mu_{v_{H,\mathbf w}}.
\end{equation}

\subsection{Comparison with the residual vector}\label{sec:bm-stop}

The purpose of the next lemma is to reduce the Brownian quantity $D_v(H)$ to an average of squared $L^2$ norms of derivatives of Gaussian smoothings. We first compare the integral containing $1/\overline M_{t(H)}^x$ with one containing $1/f_{v,e^{-H}}$, using the Gaussian kernel and earlier values of the martingale. We then apply the stopped self-similar decomposition. Rescaling each stopped piece by $e^H$ changes the smoothing scale from $h(H)=e^{-H}H^{-1/(2\tau)}$ to $H^{-1/(2\tau)}$, and the corresponding projection vector is $v_H$. The factor $e^{-2H}$ in $D_v(H)$ cancels the factors introduced by rescaling. A final weighted density comparison gives the bound by $\mathbb E[F_H(v_H)]$, where this expectation averages over the stopped coding words. The distribution of these residual vectors will be studied in the next section.

\begin{lemma}\label{lem:bm-stop}
Let $\mu$ belong to $(\mathrm I)_\tau$ or $(\mathrm{II})_\tau$,
let $1+2\tau\le s<2$, and fix a compact set
$\mathcal V\subset\mathcal D$. Use the scale $h(H)=e^{-H}H^{-1/(2\tau)}$ and the parameters
$H_0,h_0,H_\varepsilon$ from Subsection \ref{sec:stop}, with
$$
 D_v(H)=e^{-2H}\int_\R\mathbb E
 \frac{|f'_{v,h(H)}(x+B_{t(H)})|^2}
      {\overline M_{t(H)}^x}\dd x,
 \qquad t(H)=h_0^2-h(H)^2,
$$
and
$$
 F_H(v')=\|f'_{v',H^{-1/(2\tau)}}\|_2^2,\qquad v'\in\mathcal D.
$$
Here the $L^2$ norm is with respect to Lebesgue measure. There is a
constant $C<\infty$ such that, for every $v\in\mathcal V$,
$0<\varepsilon<h_0$ and $H_0\le H\le H_\varepsilon$,
\begin{equation}\label{eq:bm-stop-proof}
 D_v(H)\le C\mathbb E[F_H(v_H)].
\end{equation}
Thus
$$
 \mathcal{A}_{v,\varepsilon}
 \leq C\int_{H_0}^\infty H^{\kappa-1-1/\tau}
       \mathbb E[F_H(v_H)]\dd H.
$$
The constant depends only on the IFSs, their probability
vectors, $s$, $\tau$ and $\mathcal V$.
\end{lemma}

\begin{proof}
We first prove that, for $0\le t\le T$ and $a\ge h_t$,
\begin{equation}\label{eq:g-recip}
 \int_\R\mathbb E\frac{|f'_{v,h_t}(x+B_t)|^2}{\overline M_t^x}\dd x
 \le C_1\int_\R\frac{|f'_{v,h_t}(y)|^2}{f_{v,a}(y)}\dd y,
\end{equation}
with the universal constant $C_1=\sqrt{120}$. This is the key part of the argument. The idea is to bound the reciprocal of the running maximum $\overline M_t^x$ using the martingale at two earlier times. We choose these times so that the corresponding Gaussian smoothing scales are $2a$ and $4a$.

Let $h=h_t$. First suppose that $a\le h_0/4$, and set
$$
 t^{(1)}=h_0^2-4a^2,\qquad t^{(2)}=h_0^2-16a^2,\qquad
 X_i=\frac{B_{t^{(i)}}-B_t}{a}\quad(i=1,2).
$$
Then $0\le t^{(2)}<t^{(1)}<t$, $h_{t^{(1)}}=2a$ and $h_{t^{(2)}}=4a$.
We complete the square in the Gaussian kernel, which gives
\begin{equation}\label{eq:bm-gauss}
 f_{v,a}(y)\le k\exp\!\left(\frac{z^2}{2(k^2-1)}\right)
 f_{v,ka}(y+az),\qquad k>1.
\end{equation}
We apply \eqref{eq:bm-gauss} at $y=x+B_t$ with $(k,z)=(2,X_1)$ and $(4,X_2)$ so that:
\begin{equation}\label{eq:bm-two}
 \frac1{\overline M_t^x}
 \le\frac1{\sqrt{M_{t^{(1)}}^xM_{t^{(2)}}^x}}
 \le\frac{\sqrt8}{f_{v,a}(x+B_t)}
       \exp\!\left(\frac{X_1^2}{12}+\frac{X_2^2}{60}\right).
\end{equation}
The exponential factor has bounded expectation, uniformly in $t$ and $a$. Indeed, by independence of Brownian increments, $X_1$ and $X_2-X_1$ are independent centred Gaussians with variances $4-h^2/a^2\le4$ and $12$, respectively. We integrate first in $X_2-X_1$ and then in $X_1$ to have
\begin{equation}\label{eq:bm-exp}
 \mathbb E\exp\!\left(\frac{X_1^2}{12}+\frac{X_2^2}{60}\right)
 =\sqrt{\frac53}\,\mathbb E e^{X_1^2/9}
 =\frac{\sqrt{15}}{\sqrt{1+2h^2/a^2}}
 \le\sqrt{15}.
\end{equation}
We multiply \eqref{eq:bm-two} by $|f'_{v,h}(x+B_t)|^2$ and integrate in $x$. With $y=x+B_t$, we have
$$
 \int_\R\frac{|f'_{v,h}(x+B_t)|^2}{\overline M_t^x}\dd x
 \le\sqrt8\exp\!\left(\frac{X_1^2}{12}+\frac{X_2^2}{60}\right)
       \int_\R\frac{|f'_{v,h}(y)|^2}{f_{v,a}(y)}\dd y.
$$
Taking expectation and using \eqref{eq:bm-exp}, we have \eqref{eq:g-recip} with $C_1=\sqrt{120}$. For $a>h_0/4$, the Gaussian upper bound implies
$$
 \frac1{\overline M_t^x}\le h_0
 \le\frac4{\sqrt{2\pi}\,f_{v,a}(x+B_t)},
$$
and the substitution $y=x+B_t$ proves \eqref{eq:g-recip}.

We next show that the denominator $f_{v',1}$ in the right-hand side can be removed, uniformly over the residual vectors $v'\in\widetilde{\mathcal V}$. We use the compact support of the projected measures to bound $f_{v',1}$ from below on a fixed interval, and estimate the tails separately.

By compactness, there are constants $a<b$ depending only on the IFSs and $\mathcal V$ such that all $\mu_{v'}$ for $v'\in\widetilde{\mathcal V}$ are supported on the interval $[a,b]$. We prove, uniformly in $v'$, that
\begin{equation}\label{eq:weighted-bd}
 \int_\R\frac{|f'_{v',r}|^2}{f_{v',1}}\dd x
 \le C_2\|f'_{v',r}\|_2^2,\qquad 0<r\le1/2.
\end{equation}
Here $C_2>0$ depends only on the interval $[a,b]$.
For $x\in[a-1,b+1]$ and $y\in[a,b]$, we have $|x-y|\le b-a+1$, so
$$
 f_{v',1}(x)=\int_a^b\Gamma_1(x-y)\dd\mu_{v'}(y)
 \ge\Gamma_1(b-a+1)=c_2>0.
$$
For every $x\in\R$, by Cauchy-Schwarz,
$$
 \frac{|f'_{v',r}(x)|^2}{f_{v',1}(x)}
 =\frac{\left|\int_a^b\Gamma_r'(x-y)\dd\mu_{v'}(y)\right|^2}
        {\int_a^b\Gamma_1(x-y)\dd\mu_{v'}(y)}
 \le\int_a^b\frac{|\Gamma_r'(x-y)|^2}{\Gamma_1(x-y)}\dd\mu_{v'}(y).
$$
Since $|x-y|>1$ for $x\notin[a-1,b+1]$ and $y\in[a,b]$, we have
$$
 \begin{aligned}
 \int_{\R\setminus[a-1,b+1]}
 \frac{|f'_{v',r}(x)|^2}{f_{v',1}(x)}\dd x
 &\le\int_a^b\int_{\R\setminus[a-1,b+1]}
       \frac{|\Gamma_r'(x-y)|^2}{\Gamma_1(x-y)}\dd x\dd\mu_{v'}(y)\\
 &\le\int_{|z|>1}\frac{|\Gamma_r'(z)|^2}{\Gamma_1(z)}\dd z
 \le C_2'.
 \end{aligned}
$$
Here $C_2'>0$ is an absolute constant.

To absorb this tail bound into $\|f'_{v',r}\|_2^2$, we use that this norm is uniformly bounded below. For $x\in[b+1/2,b+1]$, we have
$$
 -f'_{v',1}(x)=\int_a^b(x-y)\Gamma_1(x-y)\dd\mu_{v'}(y)
 \ge(x-b)f_{v',1}(x)\ge\frac12 f_{v',1}(x)\ge\frac{c_2}2.
$$
Since $f'_{v',1}=f'_{v',r}*\Gamma_{\sqrt{1-r^2}}$, by Jensen's inequality we have
$$
 \frac{c_2^2}{8}\le\|f'_{v',1}\|_2^2\le\|f'_{v',r}\|_2^2.
$$
Therefore,
$$
 \begin{aligned}
 \int_\R\frac{|f'_{v',r}(x)|^2}{f_{v',1}(x)}\dd x
 &\le c_2^{-1}\int_{a-1}^{b+1}|f'_{v',r}(x)|^2\dd x+C_2'\\
 &\le\left(c_2^{-1}+\frac{8C_2'}{c_2^2}\right)\|f'_{v',r}\|_2^2,
 \end{aligned}
$$
which proves \eqref{eq:weighted-bd} with $C_2=c_2^{-1}+8C_2'/c_2^2$.

We now apply these two estimates to the stopping decomposition \eqref{eq:st-mix}. The first estimate replaces the running maximum by the smoothed density at scale $e^{-H}$. We then use Cauchy-Schwarz to bound the resulting quotient by the average over the stopped cylinders.

Fix $H\in[H_0,H_\varepsilon]$.
Since $h(H)\le e^{-H}$, \eqref{eq:g-recip} implies
$$
 \int_\R\mathbb E
 \frac{|f'_{v,h(H)}(x+B_{t(H)})|^2}{\overline M_{t(H)}^x}\dd x
 \le C_1\int_\R\frac{|f'_{v,h(H)}(y)|^2}{f_{v,e^{-H}}(y)}\dd y.
$$
Let $z_{\mathbf w}=e^H(y-t_{v,\mathbf w})$.
By \eqref{eq:st-mix}, $h(H)=e^{-H}H^{-1/(2\tau)}$ and Gaussian scaling,
$$
 e^{-H} f'_{v,h(H)}(y)
 =e^H\sum_{\mathbf w\in\mathcal P_H}
 m_{\mathbf w}f'_{v_{H,\mathbf w},H^{-1/(2\tau)}}(z_{\mathbf w}),
$$
and
$$
 f_{v,e^{-H}}(y)=e^H\sum_{\mathbf w\in\mathcal P_H}
 m_{\mathbf w}f_{v_{H,\mathbf w},1}(z_{\mathbf w}).
$$
By these identities and Cauchy-Schwarz, we have
$$
 \begin{aligned}
 \int_\R\frac{|e^{-H} f'_{v,h(H)}(y)|^2}{f_{v,e^{-H}}(y)}\dd y
 &=e^H\int_\R
 \frac{\left|\sum_{\mathbf w\in\mathcal P_H}m_{\mathbf w}
 f'_{v_{H,\mathbf w},H^{-1/(2\tau)}}(z_{\mathbf w})\right|^2}
 {\sum_{\mathbf w\in\mathcal P_H}m_{\mathbf w}
 f_{v_{H,\mathbf w},1}(z_{\mathbf w})}\dd y\\
 &\le e^H\sum_{\mathbf w\in\mathcal P_H}m_{\mathbf w}
 \int_\R\frac{|f'_{v_{H,\mathbf w},H^{-1/(2\tau)}}(z_{\mathbf w})|^2}
 {f_{v_{H,\mathbf w},1}(z_{\mathbf w})}\dd y\\
 &=\sum_{\mathbf w\in\mathcal P_H}m_{\mathbf w}
 \int_\R\frac{|f'_{v_{H,\mathbf w},H^{-1/(2\tau)}}(z)|^2}
 {f_{v_{H,\mathbf w},1}(z)}\dd z.
 \end{aligned}
$$
In the last equality, we use $y=e^{-H}z+t_{v,\mathbf w}$, so $e^H\dd y=\dd z$. Since $\tau<1/2$ and $H\ge H_0\ge e$, we have $H^{-1/(2\tau)}\le e^{-1}<1/2$. By \eqref{eq:weighted-bd}, applied to $v_{H,\mathbf w}\in\widetilde{\mathcal V}$, we have $D_v(H)\le C_3\mathbb E[F_H(v_H)]$, where $C_3=C_1C_2$.
This proves \eqref{eq:bm-stop-proof}.
By \eqref{eq:bm-scale-proof}, we then have for $C_3'=C_\kappa'C_3$:
\begin{equation}\label{eq:st-scale}
 \mathcal{A}_{v,\varepsilon}
 \le C_3'\int_{H_0}^\infty H^{\kappa-1-1/\tau}
       \mathbb E[F_H(v_H)]\dd H.
\end{equation}
\end{proof}

\section{Random stopping scales and the energy estimate}\label{sec:sob}

The earlier section reduces the logarithmic moment bound to an average over random residue  vectors. The task here is to control this average using the Diophantine assumption, where we prove a non-concentration estimate. The planar and product cases require different arguments for the interval estimate, but the scale averaging and energy estimates apply to both.

\begin{proposition}\label{prop:bm-ren}
Let $\mu$ belong to $(\mathrm{I})_\tau$ or $(\mathrm{II})_\tau$,
and suppose that $I_s(\mu)<\infty$ for some $1+2\tau\le s<2$.
Fix a compact set $\mathcal V\subset\mathcal D$.
Write $\kappa=(s-1)/(2\tau)\ge1$.
There are $h_0\in(0,1)$ and $C<\infty$ such that
$$
 \sup_{v\in\mathcal V}\sup_{0<h<h_0}
 \int_\R f_{v,h}\bigl(\log^+(h_0f_{v,h})\bigr)^\kappa\dd x
 \leq C\bigl(1+I_s(\mu)\bigr).
$$
The constants depend only on the IFSs, their probability
vectors, $s$, $\tau$ and $\mathcal V$.
\end{proposition}

To prove this, take $h_0=h(H_0)$ with $H_0$ sufficiently large.
We first use the random stopping
scale $Y_H$ to prove \eqref{eq:bm-avg-proof}. We then prove the directional interval bound
\eqref{eq:dir-mass-proof} separately in the two cases, using the angular
walk in $(\mathrm I)$ and a renewal sum in $(\mathrm{II})$.
Finally, we absorb the bounded radial
factor into the smoothing scale and integrate the Fourier estimate in
$H$ to prove \eqref{eq:g-trace-int-proof}. Together with
Lemma \ref{lem:bm-ent}, these estimates prove the proposition.

\subsection{Lazy walks and the stopping scale}\label{sec:ren-loc}

We use the renewal processes and crossing counts from Subsection \ref{sec:stop}. We introduce a \emph{lazy walk}, following Borda's approach \cite[Theorem 6]{Bor21}, which at each step either stays at its current position or takes a step of the original renewal walk, each with probability $1/2$. Set $b_\ast=\max\{1,2\beta_+\}$. Independently of $\bm\omega$, let $(\mathsf B_j)_{j\ge1}$ be independent Bernoulli variables with parameter $1/2$. Define
\begin{equation}\label{eq:lazy-scale}
 \mathsf K_n=\sum_{j=1}^n\mathsf B_j,\qquad
 n_H=\lfloor H/b_\ast\rfloor,\qquad
 \mathsf N_H=\mathsf K_{n_H},\qquad
 \widetilde W_n=W_{\mathsf K_n}.
\end{equation}
Thus $\widetilde W$ is a \emph{lazy renewal walk} with independent increments having law $\frac12\delta_0+\frac12\sum_i p_i\delta_{\beta_i}$. In $(\mathrm{II})$, let $(\mathsf B'_j)_{j\ge1}$ be another independent Bernoulli sequence with parameter $1/2$, independent of $\bm\omega$ and $(\mathsf B_j)_{j\ge1}$, and define $\mathsf K'_n$ and $\widetilde W'_n=W'_{\mathsf K'_n}$ in the same way.

We use the lazy walk to randomise the stopping scale. Let
$$
Y_H=W_{\mathsf N_H+1},\qquad \widetilde v_H=v_{Y_H}.
$$
Recall that $N_t=\min\{n\ge1:W_n\ge t\}$ for $t>0$. Since the increments of $W$ are strictly positive, $N_{Y_H}=\mathsf N_H+1$. In $(\mathrm{II})$, the second walk $W'$ stops at $N'_{Y_H}=\min\{n\ge1:W'_n\ge Y_H\}$, since $W_{N_{Y_H}}=Y_H$. By \eqref{eq:ang-comp} and \eqref{eq:prod-comp}, $\widetilde v_H\in\widetilde{\mathcal V}$ for $v\in\mathcal V$. For $\theta\in\T$, let $e_\theta=(\cos 2\pi\theta,\sin 2\pi\theta)$. Write
$$
\widetilde v_H=\rho_H e_{\vartheta_H},\qquad \rho_H=|\widetilde v_H|,
$$
and let $\sigma_{v,H}$ denote the law of $\vartheta_H\in\T$.

We next compare the average over the random stopping scale $Y_H$ with an integral over the original stopping scales. This will allow us to use $\widetilde v_H$ in the scale integral from Lemma \ref{lem:bm-stop}.

For integers $m,n\ge0$, write
\begin{equation}\label{eq:bin-wt}
 c(m,n)=2^{-m}\binom mn, \quad 0\le n\le m,\qquad
 c(m,n)=0,\quad n>m.
\end{equation}
We write $\mathbb E_{\bm\omega}$ for expectation over the coding sequences only. Since $\mathsf N_H$ is independent of these sequences and $\mathbb P(\mathsf N_H=n)=c(n_H,n)$, we have
$$
\mathbb E[F_H(\widetilde v_H)]
=\sum_{n=0}^{n_H}c(n_H,n)\,
\mathbb E_{\bm\omega}[F_H(v_{W_{n+1}})].
$$
To compare this with an integral over the stopping scale, we first show that $F_H(v_t)$ does not change by more than a fixed factor within each renewal interval $[W_n,W_{n+1}]$.

For $W_n<t\le W_{n+1}$, we have $v_t=e^{t-W_{n+1}}v_{W_{n+1}}$ and $0\le W_{n+1}-t\le\beta_+$. Let $a=e^{t-W_{n+1}}$, $h=H^{-1/(2\tau)}$ and $v'=v_{W_{n+1}}$. By Gaussian scaling and the substitution $x=ay$, we have
$$
f'_{av',h}(x)=a^{-2}f'_{v',h/a}(x/a),\qquad
\|f'_{av',h}\|_2^2=a^{-3}\|f'_{v',h/a}\|_2^2.
$$
Since $h/a\ge h$, by the Gaussian semigroup identity and $L^2$ contraction under convolution, we have
$$
F_H(v_t)\le e^{3(W_{n+1}-t)}F_H(v_{W_{n+1}})
\le e^{3\beta_+}F_H(v_{W_{n+1}}).
$$
Thus we can bound the integral over each renewal interval by its value at the right endpoint, up to a uniform constant. We integrate over each renewal interval, divide by its length $W_{n+1}-W_n=\beta_{\omega_{n+1}}$, and multiply by $c(n_H,n)$. Since $n=N_t-1$ for $W_n<t<W_{n+1}$, by Fubini and the law of total expectation, we have
\begin{equation}\label{eq:lazy-avg}
 \begin{split}
 \int_0^\infty\mathbb E_{\bm\omega}\left[
   \frac{c(n_H,N_t-1)}{\beta_{\omega_{N_t}}}F_H(v_t)\right]\dd t
 &=\sum_{n=0}^{n_H}c(n_H,n)\mathbb E_{\bm\omega}\left[
   \frac1{\beta_{\omega_{n+1}}}
   \int_{W_n}^{W_{n+1}}F_H(v_t)\dd t\right]\\
 &\le e^{3\beta_+}\sum_{n=0}^{n_H}c(n_H,n)
       \mathbb E_{\bm\omega}[F_H(v_{W_{n+1}})]\\
 &=e^{3\beta_+}\mathbb E[F_H(\widetilde v_H)].
 \end{split}
\end{equation}

\subsection{Continuous averaging in scale}\label{sec:bm-ren}

The previous comparison leaves us with a scale integral involving the stopped vector $v_H$. The directional estimates below will instead apply to $\widetilde v_H=v_{Y_H}$, where $Y_H=W_{\mathsf N_H+1}$ and $\mathsf N_H$ counts the moves of the lazy walk after $\lfloor H/b_\ast\rfloor$ steps. We therefore need to replace $v_H$ by $\widetilde v_H$ in this integral. 

The next lemma shows that using this replacement we do not lose much: we lose only a uniform constant after integration in $H$. The reason is that, as $H$ varies, the binomial weights assign a uniformly positive total weight to each renewal interval at comparable scales. The Gaussian semigroup identity then compares the corresponding smoothing norms. This allows us to use the distribution of $\widetilde v_H$ to bound $\mathcal A_{v,\varepsilon}$.

\begin{lemma}\label{lem:bm-avg}
Let $\mu$ belong to $(\mathrm I)_\tau$ or $(\mathrm{II})_\tau$,
let $1+2\tau\le s<2$, and fix a compact
$\mathcal V\subset\mathcal D$.
Let $\kappa=(s-1)/(2\tau)$ and recall that
$$
 F_H(v')=\|f'_{v',H^{-1/(2\tau)}}\|_2^2,\qquad v'\in\mathcal D.
$$
For every sufficiently large $H_0$ and every $v\in\mathcal V$,
\begin{equation}\label{eq:st-occ}
 \int_{H_0}^\infty H^{\kappa-1-1/\tau}
       \mathbb E[F_H(v_H)]\dd H
 \le C\int_{H_0}^\infty H^{\kappa-1-1/\tau}
       \mathbb E[F_H(\widetilde v_H)]\dd H.
\end{equation}
With $h_0=h(H_0)$ and $0<\varepsilon<h_0$, we also prove
\begin{equation}\label{eq:bm-avg-proof}
 \mathcal A_{v,\varepsilon}\le
 C_4\int_{H_0}^\infty H^{\kappa-1-1/\tau}
 \mathbb E[F_H(\widetilde v_H)]\dd H.
\end{equation} The constant $C$ is independent of $v$ and
$\varepsilon$.
\end{lemma}

\begin{proof}
We want to compare the integral involving $F_t(v_t)$ with the average involving $F_H(\widetilde v_H)$. The idea is to fix a renewal interval containing $t$ and find a range of $H$ for which the lazy walk gives this interval a uniformly positive total weight. We then use the Gaussian semigroup to compare $F_t$ and $F_H$.

For $m\ge n\ge0$, $c(m,n)/2$ is the probability that the $(n+1)$st success occurs at time $m+1$ in the Bernoulli sequence. Thus $\sum_{m=n}^\infty c(m,n)=2$ and
\begin{equation}\label{eq:bin-occ}
 \sum_{m=n}^{2n}c(m,n)
 =2\mathbb P(\mathsf K_{2n+1}\ge n+1)=1.
\end{equation}
Thus the weights $c(m,n)$ have total mass $1$ already for $n\le m\le2n$. We use this to choose a range of $H$ comparable to $t$.

Let $\beta_-=\min_i\beta_i>0$, and recall that $\beta_+=\max_i\beta_i$ and $b_\ast\ge\max\{1,2\beta_+\}$. Take $H_0>\beta_+$, fix a coding path and let $t\ge H_0$ lie outside its renewal heights. For $n=N_t-1\ge1$, we have $\beta_-n<t<\beta_+(n+1)$. Thus $H>t$ for every $H\in[b_\ast n,b_\ast(2n+1))$. By the definition of $n_H$ and \eqref{eq:bin-occ}, we have
$$
 \int_{b_\ast n}^{b_\ast(2n+1)}c(n_H,n)\dd H
 =b_\ast\sum_{m=n}^{2n}c(m,n)=b_\ast.
$$
We next compare the smoothing norms at scales $t$ and $H$. Since $H>t$, the smoothing scale $H^{-1/(2\tau)}$ is smaller than $t^{-1/(2\tau)}$, and hence $F_H$ is at least $F_t$. Indeed, by the Gaussian semigroup identity $\Gamma_h*\Gamma_a=\Gamma_{\sqrt{h^2+a^2}}$, we have
$$
 F_t(v')
 =\bigl\|f'_{v',H^{-1/(2\tau)}}*
   \Gamma_{\sqrt{t^{-1/\tau}-H^{-1/\tau}}}\bigr\|_2^2
 \le F_H(v').
$$
We also need to compare the weights in the two scale integrals. On $[b_\ast n,b_\ast(2n+1))$, we have $H\le(3b_\ast/\beta_-)t$ and $\kappa-1-1/\tau=(s-3)/(2\tau)-1<0$. Since $\beta_{\omega_{N_t}}\le\beta_+$ and $F_H(v_t)\ge F_t(v_t)$, the preceding bounds give
$$
 \int_{H_0}^\infty H^{\kappa-1-1/\tau}
 \frac{c(n_H,N_t-1)}{\beta_{\omega_{N_t}}}F_H(v_t)\dd H
 \ge c_1\,t^{\kappa-1-1/\tau}F_t(v_t).
$$
Here $c_1=\frac{b_\ast}{\beta_+}\left(\frac{3b_\ast}{\beta_-}\right)^{\kappa-1-1/\tau}>0$. This constant is independent of $t$, $v$ and $\varepsilon$.

We have therefore bounded the original integrand at each $t$ by an integral over $H$. We integrate this inequality over $t\ge H_0$ and take $\mathbb E_{\bm\omega}$. We extend the $t$-integral on the left to $(0,\infty)$. By Fubini and \eqref{eq:lazy-avg}, we have \eqref{eq:st-occ} with $C=e^{3\beta_+}/c_1$. By \eqref{eq:st-occ} and \eqref{eq:st-scale}, we have \eqref{eq:bm-avg-proof}.
\end{proof}

\subsection{Non-concentration of the angles $\widetilde v_H$}\label{sec:res-law}

We now use the Diophantine assumptions to control the distribution of the angle of $\widetilde v_H$. Writing $\delta_H=H^{-1/(2\tau)}$, our aim is to show that every arc of length $\delta_H$ has probability at most $C\delta_H$, uniformly over the initial projection vectors $v\in\mathcal V$. This is the same scale as the Gaussian smoothing in $F_H$. The bound therefore gives the control of directional concentration needed to compare $\mathbb E[F_H(\widetilde v_H)]$ with an angular integral in the energy estimate.

We use the Fourier transform
$\widehat\lambda(\xi)=\int e^{-2\pi i\xi\cdot x}\dd\lambda(x)$
for a finite measure $\lambda$ on $\R^d$, and the same convention
for integrable functions. In particular,
$\widehat\Gamma_h(\xi)=e^{-2\pi^2h^2\xi^2}$. The following estimate is similar to one in
Borda \cite[proof of Theorem 6]{Bor21}.

\begin{lemma}\label{lem:exp-sum}
Let $\tau>0$ and let $\bm\gamma=(\gamma_j)_{j=1}^r\in\R^r$, $r\ge1$,
satisfy $(D)_\tau$ with constant $c_{\bm\gamma}>0$.
Write $d_{\bm\gamma}(t)=\max_j\|t\gamma_j\|_{\T}$.
For every $\eta>0$, there is $C_\eta>0$ such that
\begin{equation}\label{eq:borda-sum}
 \sum_{|k|\le N}e^{-\eta n d_{\bm\gamma}(k)^2}
 \le C_\eta(1+Nn^{-1/(2\tau)}),\qquad n\ge1,\ N\ge0.
\end{equation}
The constant depends only on $\eta$, $\tau$, $\bm\gamma$ and
$c_{\bm\gamma}$.
\end{lemma}

\begin{proof}
Fix $\eta>0$.
For $n\ge1$, $N\ge0$ and $R>0$, if distinct integers $k,\ell$ satisfy
$\sqrt n\,d_{\bm\gamma}(k)\le R$ and $\sqrt n\,d_{\bm\gamma}(\ell)\le R$, then
$$
 c_{\bm\gamma}|k-\ell|^{-\tau}
 \le d_{\bm\gamma}(k-\ell)\le2R/\sqrt n.
$$
Thus these integers are at least
$(c_{\bm\gamma}\sqrt n/(2R))^{1/\tau}$ apart, so
$$
 \#\{k\in\Z:|k|\le N,\ \sqrt n\,d_{\bm\gamma}(k)\le R\}
 \le1+C_1Nn^{-1/(2\tau)}R^{1/\tau},\qquad
 C_1=2(2/c_{\bm\gamma})^{1/\tau}.
$$
Since $e^{-\eta x^2}=\int_x^\infty2\eta R e^{-\eta R^2}\dd R$
for $x\ge0$, we have
$$
 \sum_{|k|\le N}e^{-\eta n d_{\bm\gamma}(k)^2}
 \le2\eta\int_0^\infty R e^{-\eta R^2}
       (1+C_1Nn^{-1/(2\tau)}R^{1/\tau})\dd R
 \le C_\eta(1+Nn^{-1/(2\tau)}).
$$
Here
$C_\eta=1+2\eta C_1\int_0^\infty R^{1+1/\tau}e^{-\eta R^2}\dd R<\infty$.
It depends only on $\eta$, $\tau$ and $c_{\bm\gamma}$.
\end{proof}

We prove the directional bound separately in the two settings. In $(\mathrm I)$, the angle is given by a lazy random walk on $\T$. Lemma \ref{lem:exp-sum} bounds partial sums of the absolute values of its Fourier coefficients, and Gaussian smoothing converts this into the required interval estimate.

In $(\mathrm{II})$, the logarithm of the absolute slope is determined by the difference between the two logarithmic contraction walks when the second crosses the level reached by the first. We first estimate the probability that the difference of the two lazy walks is in a short interval at fixed step counts. We also need decay away from its mean so that this estimate can be summed over the step count of the second walk. This leads to a renewal bound that controls the difference at the crossing time. Finally, since the residual directions stay away from the coordinate axes, we transfer the interval bound for logarithmic absolute slopes to an interval bound for angles.

\begin{proposition}[Non-concentration]\label{prop:st-en}
Let $\mu$ belong to $(\mathrm I)_\tau$ or $(\mathrm{II})_\tau$, and
let $\sigma_{v,H}$ be the law of the angle $\vartheta_H\in\T$ defined by
$$
 \widetilde v_H=\rho_H e_{\vartheta_H},\qquad
 \rho_H=|\widetilde v_H|,\qquad
 e_\theta=(\cos 2\pi\theta,\sin 2\pi\theta).
$$
For every compact $\mathcal V\subset\mathcal D$, there are
$H_0>1$ and $C_0>0$ such that
\begin{equation}\label{eq:dir-mass-proof}
 \sigma_{v,H}(J)\le C_0H^{-1/(2\tau)}
\end{equation}
for $H\ge H_0$, $v\in\mathcal V$ and every arc $J\subset\T$
with $|J|=H^{-1/(2\tau)}$.
The constants depend only on $\tau$, $\mathcal V$, the IFSs and
their probability vectors.
\end{proposition}

\begin{proof}
\textbf{Case $(\mathrm I)$.} Fix $v\in\mathcal V$ and take $H$ sufficiently large that $H/(2b_\ast)\le n_H\le H$ and $n_H>1$. Let
$$
\lambda=\sum_i p_i\delta_{\alpha_i},\qquad \lambda_0=\tfrac12(\delta_0+\lambda)
$$
on $\T$, and write $v=|v|e_{\theta_v}$. Since $N_{Y_H}=\mathsf N_H+1$ and $\mathsf N_H=\mathsf K_{n_H}$, we have
$$
\operatorname{Law}(\Theta_{N_{Y_H}})=\lambda_0^{*n_H}*\lambda,\qquad \sigma_{v,H}=\operatorname{Law}(\theta_v-\Theta_{N_{Y_H}}).
$$
Thus the angle law is a translate and reflection of a lazy random walk on $\T$. We use the Diophantine condition to bound its Fourier coefficients.

For $k\in\Z$, we have $1-\operatorname{Re}\widehat\lambda(k)\ge8\min_i p_i\,d_{\bm\alpha}(k)^2$, and hence
$$
|\widehat\lambda_0(k)|^2\le\tfrac12(1+\operatorname{Re}\widehat\lambda(k))\le1-4\min_i p_i\,d_{\bm\alpha}(k)^2.
$$
Thus, with $c_1=2\min_i p_i>0$,
$$
|\widehat\sigma_{v,H}(k)|=|\widehat\lambda_0(k)|^{n_H}|\widehat\lambda(k)|\le e^{-c_1n_Hd_{\bm\alpha}(k)^2}.
$$

By Lemma \ref{lem:exp-sum} with $\bm\gamma=\bm\alpha$, $\eta=c_1$ and $n=n_H$, we have
\begin{equation}\label{eq:ang-sum}
\sum_{|k|\le N}|\widehat\sigma_{v,H}(k)|\le C_1(1+Nn_H^{-1/(2\tau)}),\qquad N\ge0,
\end{equation}
where $C_1=C_{c_1}$ is the constant from that lemma for $\bm\alpha$.

We now use this Fourier estimate to bound the mass of arcs of length $H^{-1/(2\tau)}$, that is, we smooth the angle law at scale $n_H^{-1/(2\tau)}$ and show that the density is uniformly bounded.

For $\delta>0$, let $\Gamma_\delta^\T(t)=\sum_{\ell\in\Z}\Gamma_\delta(t+\ell)$ on $\T$, whose Fourier coefficients are $\widehat{\Gamma}_\delta^\T(k)=e^{-2\pi^2\delta^2k^2}$. We split the frequencies into $|k|\le n_H^{1/(2\tau)}$ and $2^jn_H^{1/(2\tau)}<|k|\le2^{j+1}n_H^{1/(2\tau)}$, $j\ge0$. By \eqref{eq:ang-sum}, we have
\[
 \begin{split}
 \|\sigma_{v,H}*\Gamma_{n_H^{-1/(2\tau)}}^\T\|_\infty
 &\le\sum_{k\in\Z}e^{-2\pi^2n_H^{-1/\tau}k^2}
              |\widehat\sigma_{v,H}(k)|\\
 &\le C_1\left(2+\sum_{j=0}^\infty
            (1+2^{j+1})e^{-2\pi^2 4^j}\right)
 =:C_1'<\infty.
 \end{split}
\]

Let $J$ be an arc of length $H^{-1/(2\tau)}$ and centre $x$. Since $n_H\le H$, for $t\in J$ we have
$$
\Gamma_{n_H^{-1/(2\tau)}}^\T(x-t)\ge\frac{e^{-1/8}}{\sqrt{2\pi}\,n_H^{-1/(2\tau)}}.
$$
Thus the bound for the smoothed density gives
$$
\sigma_{v,H}(J)\le\sqrt{2\pi}e^{1/8}C_1'n_H^{-1/(2\tau)}\le C_1''H^{-1/(2\tau)},
$$
where $C_1''=\sqrt{2\pi}e^{1/8}(2b_\ast)^{1/(2\tau)}C_1'$. This proves the proposition in $(\mathrm I)$ with $C_0=C_1''$.

\textbf{Case $(\mathrm{II})$.} Recall that the direction of $\widetilde v_H$ is determined by the overshoot $\Delta_H$ of the second contraction walk. We first estimate the difference $D_{n,\ell}=\widetilde W_n-\widetilde W'_\ell$ at fixed step counts. We then sum over $\ell$ to obtain a renewal estimate, which will be used to control $\Delta_H$ and hence the direction of $\widetilde v_H$.

We begin by proving a real-frequency estimate from \eqref{eq:borda-sum}. If $\bm\gamma$ satisfies $(D)'_\tau$ with constant $c_{\bm\gamma}>0$, then
\begin{equation}\label{eq:dio-small}
d_{\bm\gamma}(t)\ge 2^{-(\tau+1)}c_{\bm\gamma}\min\{|t|,|t|^{-\tau}\},\qquad t\ne0.
\end{equation}
Indeed, for $0<|t|<1$, let $k=\lceil1/|t|\rceil$. Since $1\le k|t|<2$ and $k\le2/|t|$, we have $k d_{\bm\gamma}(t)\ge d_{\bm\gamma}(kt)\ge2^{-\tau}c_{\bm\gamma}$. For $|t|\ge1$, use $(D)'_\tau$ and symmetry.

Let $\bm\beta$ be the joint vector of logarithmic contractions from the two IFSs, and let $b=\beta_1>0$. By \eqref{eq:dio-small}, the vector $\bm\beta/b$ satisfies $(D)_\tau$. Let $X_j,X'_j$ be the increments of $\widetilde W,\widetilde W'$, and write $D_{n,\ell}=\widetilde W_n-\widetilde W'_\ell$ and $a_1=\mathbb E X_1>0$, $a_2=\mathbb E X'_1>0$. Choose
$$
C_\ast>\max\{2,4\beta'_+/a_1,4\beta_+/a_2\}.
$$
For $\varphi_1(t)=\mathbb E e^{-2\pi itX_1}$ and $\varphi_2(t)=\mathbb E e^{-2\pi itX'_1}$, the lazy increment laws imply
$$
|\varphi_1(t)|^n|\varphi_2(t)|^\ell\le e^{-3\eta_\ast n d_{\bm\beta}(t)^2},\qquad C_\ast^{-1}n\le\ell\le C_\ast n,
$$
where $\eta_\ast>0$ depends only on the increment laws.

We first reduce the real-frequency integral to a sum over integers, so that we can apply \eqref{eq:borda-sum}. We use the coordinate $b=\beta_1$ to control the distance from $t$ to the nearest $k/b$, and the Diophantine condition to control the values at these points. Partition first $\R$ into intervals
$$
I_k=[(k-\tfrac12)/b,(k+\tfrac12)/b),\qquad k\in\Z.
$$
The constants introduced below are fixed and depend only on the IFSs,
their probability vectors, $\tau$ and $\mathcal V$.
Set $C_{\bm\beta}=1+\max\{\beta_+,\beta'_+\}/b$. For $t\in I_k$, the triangle inequality on $\T$ implies
$$
b|t-k/b|=\|tb\|_{\T}\le d_{\bm\beta}(t),\qquad d_{\bm\beta}(k/b)\le C_{\bm\beta}d_{\bm\beta}(t).
$$
Set $\eta=\eta_\ast/(2C_{\bm\beta}^2)$. For $t\in I_k$, we have
$$
d_{\bm\beta}(t)^2\ge\frac{b^2}{2}(t-k/b)^2+\frac{d_{\bm\beta/b}(k)^2}{2C_{\bm\beta}^2}.
$$
Hence
\[
 \begin{split}
 \int_{I_k}e^{-\eta_\ast n d_{\bm\beta}(t)^2}\dd t
 &\le e^{-\eta n d_{\bm\beta/b}(k)^2}
       \int_{\R}e^{-\eta_\ast n b^2(t-k/b)^2/2}\dd t\\
 &=\frac{\sqrt{2\pi}}{b\sqrt{\eta_\ast}}\,
       n^{-1/2}e^{-\eta n d_{\bm\beta/b}(k)^2}\\
 &=C_2n^{-1/2}e^{-\eta n d_{\bm\beta/b}(k)^2},
 \end{split}
\]
where $C_2=\sqrt{2\pi}/(b\sqrt{\eta_\ast})$.

Let $C_\eta$ be the constant in \eqref{eq:borda-sum} for $\bm\beta/b$, and set
$$
\varepsilon_1=\min\{1/2,1/(2b)\},\qquad T_n=\varepsilon_1 n^{1/(2\tau)}.
$$
For all sufficiently large $n$, we have $bT_n+1\le n^{1/(2\tau)}$. Since $d_{\bm\beta}(k/b)=d_{\bm\beta/b}(k)$, by \eqref{eq:borda-sum} we have
\begin{equation}\label{eq:real-fourier}
\int_{-T_n}^{T_n} e^{-\eta_\ast n d_{\bm\beta}(t)^2}\dd t\le C_2n^{-1/2}\sum_{|k|\le bT_n+1}e^{-\eta n d_{\bm\beta/b}(k)^2}\le C_2'n^{-1/2}.
\end{equation}
Here $C_2'=2C_2C_\eta$.

We now use \eqref{eq:real-fourier} to estimate the probability that $D_{n,\ell}$ lies in a short interval. We also need decay when the interval is far from the mean $na_1-\ell a_2$, in order to sum over $\ell$. For this, we split the centred difference into two independent blocks, using the concentration estimate for one block and a large-deviation estimate for the other.

Fix $n,\ell$ with $C_\ast^{-1}n\le\ell\le C_\ast n$, and let $J$ be an interval of length $n^{-1/(2\tau)}$ with centre $x$. Define
$$
Q_{1,n,\ell}=\sum_{r=1}^{\lfloor n/2\rfloor}(X_r-a_1)-\sum_{r=1}^{\lfloor\ell/2\rfloor}(X'_r-a_2),
$$
and
$$
Q_{2,n,\ell}=\sum_{r=\lfloor n/2\rfloor+1}^{n}(X_r-a_1)-\sum_{r=\lfloor\ell/2\rfloor+1}^{\ell}(X'_r-a_2).
$$
Since the two walks are independent and have independent increments, $Q_{1,n,\ell}$ and $Q_{2,n,\ell}$ are independent. We have $D_{n,\ell}=na_1-\ell a_2+Q_{1,n,\ell}+Q_{2,n,\ell}$. For sufficiently large $n$, we have $\lfloor n/2\rfloor\ge n/3$ and $\lfloor\ell/2\rfloor\ge\ell/3$. Hence, for $j=1,2$ and every $t\in\R$,
$$
|\mathbb E e^{-2\pi itQ_{j,n,\ell}}|\le\bigl(|\varphi_1(t)|^n|\varphi_2(t)|^\ell\bigr)^{1/3}\le e^{-\eta_\ast n d_{\bm\beta}(t)^2}.
$$
Since $T_n n^{-1/(2\tau)}=\varepsilon_1<1$, by Esseen's concentration inequality \cite{Ess66} at scale $T_n^{-1}$ and \eqref{eq:real-fourier}, there is an absolute constant $C_3>0$ such that
\begin{equation}\label{eq:blk-mass}
\sup_{|I|=n^{-1/(2\tau)}}\mathbb P(Q_{j,n,\ell}\in I)\le\frac{C_3}{T_n}\int_{-T_n}^{T_n}e^{-\eta_\ast n d_{\bm\beta}(t)^2}\dd t\le C_3'n^{-1/2-1/(2\tau)},\qquad j=1,2,
\end{equation}
where $C_3'=C_3C_2'/\varepsilon_1$. By independence and conditioning on $Q_{2,n,\ell}$, we have
\begin{equation}\label{eq:lift-mass}
\mathbb P(D_{n,\ell}\in J)\le\sup_{y\in\R}\mathbb P(Q_{1,n,\ell}\in J-y)\le C_3'n^{-1/2-1/(2\tau)}.
\end{equation}

The estimate \eqref{eq:lift-mass} does not use the fact that $D_{n,\ell}$ has mean $na_1-\ell a_2$. Summing it over the $O(n)$ comparable values of $\ell$ would give a bound that is too large. We therefore use Hoeffding's inequality to gain additional decay when $x$ is far from this mean.

Let $z=x-na_1+\ell a_2$. If $|z|>\sqrt n$ and $D_{n,\ell}\in J$, then $|Q_{j,n,\ell}|\ge|z|/4$ for at least one $j$. We condition on $Q_{j,n,\ell}$ and apply \eqref{eq:blk-mass} to the other block. By Hoeffding's inequality, we have
$$
\mathbb P(D_{n,\ell}\in J)\le C_3'n^{-1/2-1/(2\tau)}\sum_{j=1}^2\mathbb P(|Q_{j,n,\ell}|\ge|z|/4)\le C_3''n^{-1/2-1/(2\tau)}e^{-c_3z^2/n},
$$
where $C_3''=4C_3'$ and $c_3=(8(\beta_+^2+C_\ast(\beta'_+)^2))^{-1}$. For $|z|\le\sqrt n$, we have $e^{-c_3z^2/n}\ge e^{-c_3}$. We combine this with \eqref{eq:lift-mass} to obtain
\begin{equation}\label{eq:lazy-local}
\mathbb P(D_{n,\ell}\in J)\le C_3'''n^{-1/2-1/(2\tau)}e^{-c_3(x-na_1+\ell a_2)^2/n}.
\end{equation}
Here $C_3'''=e^{c_3}C_3''$.

We now sum over all $\ell\ge0$. The estimate \eqref{eq:lazy-local} applies only when $n$ and $\ell$ are comparable. For the remaining indices, the two walks can be close only if at least one of them deviates substantially from its mean. We use Hoeffding's inequality to bound these contributions.

Set $R_0=\beta_++\beta'_++2$. For sufficiently large $n$, if $\ell<C_\ast^{-1}n$ and $|D_{n,\ell}|\le R_0$, then $\widetilde W_n\le\beta'_+\ell+R_0\le a_1n/2$. Similarly, if $\ell>C_\ast n$ and $|D_{n,\ell}|\le R_0$, then $\widetilde W'_\ell\le\beta_+n+R_0\le a_2\ell/2$. Let
$$
c_4=\frac{a_1^2}{2\beta_+^2},\qquad c_4'=\frac{a_2^2}{2(\beta'_+)^2}.
$$
Since $X_j\in[0,\beta_+]$ and $X'_j\in[0,\beta'_+]$, by Hoeffding's inequality we have
$$
\mathbb P(|D_{n,\ell}|\le R_0)\le\begin{cases}e^{-c_4n},&\ell<C_\ast^{-1}n,\\ e^{-c_4'\ell},&\ell>C_\ast n.\end{cases}
$$
There are at most $n$ indices in the first range. We sum the second bound as a geometric series to obtain
$$
\sum_{\ell<C_\ast^{-1}n\ \mathrm{or}\ \ell>C_\ast n}\mathbb P(|D_{n,\ell}|\le R_0)\le ne^{-c_4n}+\frac{e^{-c_4'C_\ast n}}{1-e^{-c_4'}}\le C_4e^{-\eta_1n},
$$
where
$$
\eta_1=\min\{c_4/2,C_\ast c_4'\},\qquad C_4=\frac{2}{ec_4}+\frac{1}{1-e^{-c_4'}}.
$$

For the remaining indices, we sum \eqref{eq:lazy-local}. The Gaussian factor is concentrated around $\ell=(na_1-x)/a_2$, with width of order $\sqrt n$. Thus the sum is $O(\sqrt n)$, which cancels the factor $n^{-1/2}$ in \eqref{eq:lazy-local}. To see this, we extend the sum to all $\ell\in\Z$, bound the two terms bracketing $(na_1-x)/a_2$ by $1$, and compare the two monotone tails with an integral. Thus
$$
\sum_{\ell\ge0}e^{-c_3(x-na_1+\ell a_2)^2/n}\le 2+\frac1{a_2}\int_{\R}e^{-c_3t^2/n}\dd t=2+\frac{\sqrt{\pi n/c_3}}{a_2}\le C_4'\sqrt n,
$$
where $C_4'=2+\sqrt{\pi/c_3}/a_2$. For sufficiently large $n$, we have $e^{-\eta_1n}\le n^{-1/(2\tau)}$. Thus, with $C_4''=C_3'''C_4'+C_4$, we have
\begin{equation}\label{eq:lazy-ren}
\sum_{\ell\ge0}\mathbb P(D_{n,\ell}\in J)\le C_4''n^{-1/(2\tau)}
\end{equation}
for every sufficiently large $n$ and every interval $J\subset[-R_0,R_0]$ with $|J|=n^{-1/(2\tau)}$.

We next pass from the renewal estimate \eqref{eq:lazy-ren} to $\Delta_H$. The renewal estimate concerns the lazy walks, whereas $\Delta_H$ is defined using the original second walk at its crossing time. We use the following identity to compare their occupation measures.

By conditioning on $\mathsf K'_\ell$ and summing over $\ell$, we have
\begin{equation}\label{eq:hold-id}
\sum_{\ell\ge0}\operatorname{Law}(\widetilde W'_\ell)=\sum_{r\ge0}\left(\sum_{\ell\ge r}c(\ell,r)\right)\operatorname{Law}(W'_r)=2\sum_{r\ge0}\operatorname{Law}(W'_r).
\end{equation}
These measures are locally finite, since the increments of $W'$ are bounded below by $\min_j\beta'_j>0$. Let
$$
\Delta_H=-Z'_{Y_H}=Y_H-W'_{N'_{Y_H}}\in[-\beta'_+,0].
$$

Conditional on $\mathsf N_H=r$, the second walk crosses $Y_H=W_{r+1}$ at step $r'+1$ precisely when $W'_{r'}<W_{r+1}\le W'_{r'+1}$. We sum over $r,r'$ and discard this crossing condition to prove an upper bound. Using independence of the next symbols and \eqref{eq:hold-id}, we can then apply the renewal estimate \eqref{eq:lazy-ren}. We have:
\begin{equation}\label{eq:prod-dir}
\begin{split}
\mathbb P(\Delta_H\in J)
&=\sum_{r,r'\ge0}c(n_H,r)\,\mathbb P\bigl(W'_{r'}<W_{r+1}\le W'_{r'+1},\ W_{r+1}-W'_{r'+1}\in J\bigr)\\
&\le\sum_{r,r'\ge0}c(n_H,r)\mathbb P(W_{r+1}-W'_{r'+1}\in J)\\
&=\frac12\sum_{i,j}p_ip'_j\sum_{\ell\ge0}\mathbb P\bigl(\widetilde W_{n_H}-\widetilde W'_\ell\in J-\beta_i+\beta'_j\bigr)\\
&\le C_5H^{-1/(2\tau)}
\end{split}
\end{equation}
for intervals $J$ of length $H^{-1/(2\tau)}$. Since $\Delta_H\in[-\beta'_+,0]$, we may restrict to intervals $J$ meeting $[-\beta'_+,0]$. For sufficiently large $H$, these intervals and their translates $J-\beta_i+\beta'_j$ lie in the interior of $[-R_0,R_0]$. Since $n_H\le H$, each translated interval can be enlarged to length $n_H^{-1/(2\tau)}$ inside $[-R_0,R_0]$. The last inequality therefore follows from \eqref{eq:lazy-ren} and $H/(2b_\ast)\le n_H\le H$, with $C_5=(2b_\ast)^{1/(2\tau)}C_4''$. The bound for shorter intervals follows by monotonicity.

Finally, we have the identity:
$$
\log\left|\frac{(\widetilde v_H)_1}{(\widetilde v_H)_2}\right|=\log|v_1/v_2|-\Delta_H.
$$
By \eqref{eq:prod-comp}, the directions of $\widetilde v_H$ stay uniformly away from the coordinate axes. On a compact subarc of each quadrant containing these directions, the map $\theta\mapsto\log|\cos(2\pi\theta)/\sin(2\pi\theta)|$ has derivative bounded in absolute value by a fixed $C_5'>0$. Hence an arc of length $H^{-1/(2\tau)}$ corresponds, for each sign pair, to an interval of logarithmic absolute slopes of length at most $C_5'H^{-1/(2\tau)}$. We cover it by at most $\lceil C_5'\rceil+1$ intervals of length $H^{-1/(2\tau)}$. We apply \eqref{eq:prod-dir} to these intervals and sum over the four sign pairs to obtain the claimed bound with $C_0=C_5''=4(\lceil C_5'\rceil+1)C_5$.
\end{proof}

\subsection{The energy estimate}\label{sec:g-av}

We have reduces the logarithmic moment estimate to bounding the weighted scale integral of $\mathbb E[F_H(\widetilde v_H)]$ by $I_s(\mu)$. We now prove this bound. 

Plancherel's identity allows us to write $F_H$ as an integral of $|\widehat\mu|^2$ along the projection direction. The interval bound from Proposition \ref{prop:st-en} then allows us to bound the random directional average by an integral over all angles, with a gradient term that controls the variation of $\widehat\mu$ within each short arc. Then integrating in $H$ and using polar coordinates leads to weight $|\xi|^{s-2}$ since $2\tau\kappa=s-1$. Here the compact support allows us to bound the gradient term by energy. 

\begin{lemma}\label{lem:g-trace}
Let $\mu$ belong to $(\mathrm I)_\tau$ or $(\mathrm{II})_\tau$, with
$I_s(\mu)<\infty$ for some $1+2\tau\le s<2$, and fix a compact set
$\mathcal V\subset\mathcal D$. Write $\kappa=(s-1)/(2\tau)$.
Let $H_0,C_0$ be as in Proposition \ref{prop:st-en}.
Then
\begin{equation}\label{eq:g-trace-int-proof}
 \int_{H_0}^\infty H^{\kappa-1-1/\tau}
 \mathbb E[F_H(\widetilde v_H)]\dd H
 \le C_6\int_{\R^2}|\xi|^{s-2}
 \bigl(|\widehat\mu(\xi)|^2+|\nabla\widehat\mu(\xi)|^2\bigr)\dd\xi
 \le C_6'I_s(\mu).
\end{equation}
In particular,
$$
 \int_{H_0}^\infty H^{\kappa-1-1/\tau}\mathbb E[F_H(\widetilde v_H)]\dd H\le CI_s(\mu)
$$
uniformly in $v\in\mathcal V$.
The constant depends only on $H_0$, $C_0$, $\widetilde{\mathcal V}$,
$E$, $s$ and $\tau$.
\end{lemma}

\begin{proof}
Fix $v\in\mathcal V$ and $H\ge H_0$.
By compactness of $\widetilde{\mathcal V}\subset\mathcal D$, we have
$0<\rho_-\le\rho_H\le\rho_+<\infty$ uniformly in $H,v$.
Set $C_1=4\pi^2\rho_-^{-3}$ and $c_1=4\pi^2/\rho_+^2$.
For $\rho\in[\rho_-,\rho_+]$, Plancherel's identity and the substitution $R=\rho\xi$ implies:
\[
\begin{split}
 F_H(\rho e_\theta)
 &=\frac{4\pi^2}{\rho^3}\int_\R R^2
       |\widehat\mu(R e_\theta)|^2
       e^{-4\pi^2H^{-1/\tau}R^2/\rho^2}\dd R\\
 &\le C_1\int_\R R^2e^{-c_1H^{-1/\tau}R^2}
       |\widehat\mu(R e_\theta)|^2\dd R.
\end{split}
\]

We now use the non-concentration estimate for $\sigma_{v,H}$ to replace integration over the random directions by integration over all directions. For this, we partition $\T$ into intervals of length comparable to $H^{-1/(2\tau)}$ and control the variation of $\widehat\mu$ on each interval.

Partition $\T$ into intervals $I$ of equal length $\ell$, where
$\frac12H^{-1/(2\tau)}\le\ell\le H^{-1/(2\tau)}$.
Since $|\partial_\theta\widehat\mu(R e_\theta)|
\le2\pi|R|\,|\nabla\widehat\mu(R e_\theta)|$, by the fundamental
theorem of calculus and Cauchy-Schwarz we have
$$
 \sup_{\theta\in I}|\widehat\mu(R e_\theta)|^2 \le\frac2\ell\int_I|\widehat\mu(R e_\theta)|^2\dd\theta +8\pi^2R^2\ell \int_I|\nabla\widehat\mu(R e_\theta)|^2\dd\theta.
$$
Multiplying by $\sigma_{v,H}(I)\le C_0H^{-1/(2\tau)}$ and summing
over the partition, we have
\begin{equation}\label{eq:ang-samp}
 \int_\T|\widehat\mu(R e_\theta)|^2\dd\sigma_{v,H}(\theta) \le C_2\int_\T\bigl(|\widehat\mu(R e_\theta)|^2 +H^{-1/\tau}R^2|\nabla\widehat\mu(R e_\theta)|^2\bigr) \dd\theta.
\end{equation}
Here $C_2=8\pi^2C_0$.
The bound for $F_H(\rho e_\theta)$ is uniform in $\rho\in[\rho_-,\rho_+]$.
Taking expectation with $\rho=\rho_H$ and $\theta=\vartheta_H$, and using
the law $\sigma_{v,H}$ and \eqref{eq:ang-samp}, we have
\begin{equation}\label{eq:quot-av}
 \mathbb E[F_H(\widetilde v_H)] \le C_2'\int_\R R^2e^{-c_1H^{-1/\tau}R^2} \int_\T\bigl(|\widehat\mu(R e_\theta)|^2 +H^{-1/\tau}R^2|\nabla\widehat\mu(R e_\theta)|^2\bigr) \dd\theta\dd R.
\end{equation}
Here $C_2'=C_1C_2$.

We next integrate in $H$. The choice of the weight $H^{\kappa-1-1/\tau}$ gives the power $|R|^{s-1}$ after integration, which will allow us to use the Fourier representation of $s$-energy.

For $R\ne0$, the substitution $u=H^{-1/(2\tau)}|R|$
and $2\tau\kappa=s-1$ imply
$$
\int_0^\infty H^{\kappa-1-1/\tau}R^2
       (1+H^{-1/\tau}R^2)e^{-c_1H^{-1/\tau}R^2}\dd H
 =2\tau|R|^{s-1}\int_0^\infty u^{2-s}(1+u^2)e^{-c_1u^2}\dd u=:C_3|R|^{s-1}.
$$
The integral in $u$ is finite since $1<s<2$.
We integrate \eqref{eq:quot-av} in $H$ with weight
$H^{\kappa-1-1/\tau}$ and extend the integrals on the right to $(0,\infty)$.
By the above identity and polar coordinates, we have
\begin{equation}\label{eq:quot-four}
 \int_{H_0}^\infty H^{\kappa-1-1/\tau} \mathbb E[F_H(\widetilde v_H)]\dd H \le C_3'\int_{\R^2}|\xi|^{s-2} \bigl(|\widehat\mu(\xi)|^2+|\nabla\widehat\mu(\xi)|^2\bigr)\dd\xi.
\end{equation}
Here $C_3'=C_2'C_3/\pi$. So we just need to bound the two Fourier integrals by $I_s(\mu)$. The first follows directly from the Fourier representation of energy, while the second follows by applying the same identity to the measures $x_j\mu$. By Mattila \cite[Lemma 12.12]{Ma95}:
$
 \int_{\R^2}|\xi|^{s-2}|\widehat\mu(\xi)|^2\dd\xi=C_sI_s(\mu),
$
where $C_s>0$ depends only on $s$. Let $D_E=\sup_{x\in E}|x|$. Since
$\partial_j\widehat\mu=-2\pi i\,\widehat{x_j\mu}$, the same identity,
extended to the signed measures $x_j\mu$, implies
$$
 \int_{\R^2}|\xi|^{s-2}|\nabla\widehat\mu(\xi)|^2\dd\xi =4\pi^2C_s\iint\frac{x\cdot y}{|x-y|^s}\dd\mu(x)\dd\mu(y) \leq 4\pi^2C_sD_E^2I_s(\mu).
$$
The double integral is absolutely convergent because
$|x\cdot y|\le D_E^2$ on $E\times E$ and $I_s(\mu)<\infty$.
Applying these bounds to \eqref{eq:quot-four}, we have
\eqref{eq:g-trace-int-proof} with $C_6=C_3'$ and $C_6'=C_4$, where
$C_4=C_3'C_s(1+4\pi^2D_E^2)$.
\end{proof}

\begin{proof}[Proof of Proposition \ref{prop:bm-ren}]
Choose $H_0$ sufficiently large for Lemmas \ref{lem:bm-stop},
\ref{lem:bm-avg} and \ref{lem:g-trace}, and set $h_0=h(H_0)$.
By \eqref{eq:st-scale}, Lemma \ref{lem:bm-avg} and
Lemma \ref{lem:g-trace}, we have $\mathcal{A}_{v,\varepsilon}\le C_1I_s(\mu)$
uniformly in $v\in\mathcal V$ and $0<\varepsilon<h_0$, where
$C_1>0$ depends only on the IFSs, their probability vectors, $s$,
$\tau$ and $\mathcal V$.
Since $\varepsilon$ is arbitrary, Lemma \ref{lem:bm-ent} proves the proposition.
\end{proof}

\section{Absolute continuity and dimension conservation}
\label{sec:ac}\label{sec:sets}

We now use Proposition \ref{prop:bm-ren} to complete the proofs of the measure and set theorems. The logarithmic moment bound  implies uniform integrability of the Gaussian smoothings and hence absolute continuity of every projection. In Subsection \ref{sec:pos-dim}, we prove positivity of the densities almost everywhere on the projected attractors and, under SSC, exact dimensionality of the conditional measures. Subsection \ref{sec:ssc-app} then proves set version by constructing strongly separated subsystems that approximate the original dimensions while preserving the Diophantine exponent.

\subsection{Proof of the measure theorem}\label{sec:pos-dim}

\begin{proof}[Proof of Theorem \ref{thm:ac}]
Write $\kappa=(s-1)/(2\tau)\ge1$. Fix a compact $\mathcal V\subset\mathcal D$ and choose $h_0$ as in Proposition \ref{prop:bm-ren}. By Proposition \ref{prop:bm-ren}, we have
\begin{equation}\label{eq:ren-ent-ac}
 \sup_{v\in\mathcal V}\sup_{0<h<h_0}
 \int_\R f_{v,h}\bigl(\log^+(h_0f_{v,h})\bigr)^\kappa\dd x
 \le C_0\bigl(1+I_s(\mu)\bigr)
\end{equation}
for fixed constants $C_0<\infty$ and $h_0\in(0,1)$, depending only on the IFSs, their probability vectors, $s$, $\tau$ and $\mathcal V$. We first use this estimate to pass from the Gaussian smoothings to densities of the projected measures.

Since $\log(e+z)\le C_{h_0}(1+\log^+(h_0z))$, the Gaussian smoothings satisfy
$$
 \sup_{v\in\mathcal V}\sup_{0<h<h_0}
 \int_\R f_{v,h}\bigl(\log(e+f_{v,h})\bigr)^\kappa\dd x
 \le C_{\mathcal V}
$$
for some $C_{\mathcal V}>0$. By uniform integrability and the weak convergence $f_{v,h}\dd x\rightarrow\mu_v$ as $h\downarrow0$, we have $\mu_v=f_v\dd x$. Since Gaussian approximation implies $f_{v,h}\rightarrow f_v$ almost everywhere, Fatou's lemma proves \eqref{eq:log-mom}. Since $f_v$ is supported on the compact set $E_v$, the negative part of $f_v\log f_v$ is integrable. Since $\kappa\ge1$, by Jensen's inequality and \eqref{eq:log-mom} we have
$$
 -\log\LL(E_v)\le\int_\R f_v\log f_v\dd x
 \le\int_\R f_v\log^+f_v\dd x\le C_{\mathcal V}.
$$
Thus $\inf_{v\in\mathcal V}\LL(E_v)\ge e^{-C_{\mathcal V}}>0$, and
\begin{equation}\label{eq:ent-bd}
 \sup_{v\in\mathcal V}\int_\R f_v\log^+f_v\dd x<\infty
 \quad\text{for every compact }\mathcal V\subset\mathcal D.
\end{equation}

We next prove that $f_v$ is positive almost everywhere on $E_v$. We use the stopping decomposition to find, near every point of $E_v$, a projected cylinder on which the set $\{f_v>0\}$ has positive relative Lebesgue measure.

Fix $v\in\mathcal D$. Take integer scales $H=n\ge1$. For $\mathbf w\in\mathcal P_n$, $n\ge1$, write $C_{\mathbf w}=\Psi_{\mathbf w}(E)$. By \eqref{eq:st-map},
\begin{equation}\label{eq:cyl-map}
 \pi_v\circ\Psi_{\mathbf w}
 =\Phi_{t_{v,\mathbf w},e^{-n}}\circ\pi_{v_{n,\mathbf w}}.
\end{equation}
By the stopping condition, the contraction ratios of $\Psi_{\mathbf w}$ are at most $e^{-n}$ in both cases. Hence
\begin{equation}\label{eq:cyl-sc}
 \operatorname{diam}C_{\mathbf w}\le e^{-n}\operatorname{diam}E.
\end{equation}
By \eqref{eq:ang-comp} and \eqref{eq:prod-comp}, all $v_{n,\mathbf w}$ lie in a compact set $\mathcal V_v\subset\mathcal D$ independent of $n$.

For each $v'\in\{v\}\cup\mathcal V_v$, choose a nonnegative Borel density $f_{v'}$ of $\mu_{v'}$ vanishing outside $E_{v'}$. By \eqref{eq:st-mix}, for almost every $y$,
\begin{equation}\label{eq:cyl-dec}
 f_v(y)=\sum_{\mathbf w\in\mathcal P_n}
 e^nm_{\mathbf w}
 f_{v_{n,\mathbf w}}\bigl(e^n(y-t_{v,\mathbf w})\bigr).
\end{equation}
Write $B_{v'}=\{f_{v'}>0\}$. By the Lebesgue density theorem, it suffices to show the following, as in \cite[Proposition 8.2]{JS26} and \cite[Section 7]{Ra20}:
\begin{equation}\label{eq:den-dev}
 \limsup_{r\downarrow0}\frac{\LL(B_v\cap[y-r,y+r])}{2r}>0
 \qquad(y\in E_v).
\end{equation}
By Jensen's inequality and \eqref{eq:ent-bd}, $\LL(B_{v'})\ge c_v>0$ uniformly for $v'\in\mathcal V_v$. Since $m_{\mathbf w}>0$, \eqref{eq:cyl-dec} and \eqref{eq:cyl-map} imply
$$
 \LL(B_v\cap\pi_vC_{\mathbf w})
 \ge e^{-n}\LL(B_{v_{n,\mathbf w}})\ge c_ve^{-n}.
$$
For every $y\in E_v$ and $n\ge1$, choose $\mathbf w\in\mathcal P_n$ with $y\in\pi_vC_{\mathbf w}$. By \eqref{eq:cyl-sc}, with $C_v=|v|C_2$, we have $\pi_vC_{\mathbf w}\subset[y-C_ve^{-n},y+C_ve^{-n}]$. Thus \eqref{eq:den-dev} holds, and $\mu_v\sim\LL|_{E_v}$.

We finally compute the dimensions of conditional measures under SSC. We first use the stopped cylinders to compute the dimension of $\mu$, and then compare their masses under the conditional measures.

Suppose that SSC holds. Then $\mu|_{C_{\mathbf w}}=m_{\mathbf w}(\Psi_{\mathbf w})_*\mu$. Let $X=X(\bm\omega)$ have law $\mu$ under the Bernoulli coding, with independent coordinates in $(\mathrm{II})$, and let $\mathbf w_n$ be its stopped index. By SSC and the first-crossing rule, there are $c_3,C_3>0$, depending only on the IFSs, such that, uniformly in $X,n$,
\begin{equation}\label{eq:ball-comp}
 B(X,c_3e^{-n})\cap E\subset C_{\mathbf w_n}
 \subset B(X,C_3e^{-n}).
\end{equation}
By the strong law in each coordinate, the bounds $W_{N_n}=n+O(1)$ and, in the product case, $W'_{N'_n}=n+O(1)$ from the first-crossing rule, as well as \eqref{eq:ball-comp}, $\mu$ is exact dimensional with
\begin{equation}\label{eq:src-dim}
 \frac{-\log m_{\mathbf w_n}}n\rightarrow\dim_H\mu
 \quad\text{almost surely}.
\end{equation}
By SSC, \eqref{eq:cyl-map} and disintegration, we have
$$
 (\pi_v)_*(\mu|_{C_{\mathbf w}})
 =e^nm_{\mathbf w}
   f_{v_{n,\mathbf w}}\bigl(e^n(y-t_{v,\mathbf w})\bigr)\dd y
 =f_v(y)\mu_{v,y}(C_{\mathbf w})\dd y.
$$
Since $0<f_v(y)<\infty$ for $\mu_v$-almost every $y$, we divide by $f_v(y)$ to have
\begin{equation}\label{eq:cond-mass}
 \mu_{v,y}(C_{\mathbf w})
 =\frac{e^nm_{\mathbf w}}{f_v(y)}
 f_{v_{n,\mathbf w}}\bigl(e^n(y-t_{v,\mathbf w})\bigr)
\end{equation}
for $\mu_v$-almost every $y$. Since the cylinders are countable, one may take one exceptional null set for all integer $n\ge1$. See also \cite[Lemma 7.1]{JS26} and \cite[Lemma 5.1]{Ra20}.

The formula \eqref{eq:cond-mass} expresses the conditional cylinder mass in terms of the original cylinder mass and a density of a residual projection. We need to show that the latter density contributes only a subexponential factor.

Let
$$
 Q_n=f_{v_{n,\mathbf w_n}}
 \bigl(e^n(\pi_vX-t_{v,\mathbf w_n})\bigr).
$$
Conditional on $\mathbf w_n=\mathbf w$, the variable $\Psi_{\mathbf w}^{-1}X$ has law $\mu$. By \eqref{eq:cyl-map}, $e^n(\pi_vX-t_{v,\mathbf w})$ has density $f_{v_{n,\mathbf w}}$. Hence $0<Q_n<\infty$ almost surely. At $y=\pi_vX$, \eqref{eq:cond-mass} implies
$$
 -\log\mu_{v,\pi_vX}(C_{\mathbf w_n})
 =-\log m_{\mathbf w_n}-n-\log Q_n+\log f_v(\pi_vX).
$$
Here $0<f_v(\pi_vX)<\infty$ almost surely. All $f_{v'}$, $v'\in\mathcal V_v$, are supported on a common bounded interval $I$. By \eqref{eq:ent-bd} and $z|\log z|\le1/e$ for $0<z\le1$,
$$
 \sup_n\mathbb E|\log Q_n|
 \le\sup_{v'\in\mathcal V_v}\int_\R f_{v'}|\log f_{v'}|\dd x
 \le C_{\mathcal V_v}+\frac{\LL(I)}e=:C_3'<\infty.
$$

We first prove that $\log Q_n=o(n)$ along the subsequence $n_j=j^2$, for which the probabilities obtained from Markov's inequality are summable. We then use the monotonicity of cylinder masses to pass to all $n$.

Set $n_j=j^2$. For every $\varepsilon>0$, by Markov's inequality we have
$$
 \sum_j\mathbb P(|\log Q_{n_j}|>\varepsilon n_j)
 \le\frac{C_3'}\varepsilon\sum_jn_j^{-1}<\infty.
$$
By the Borel-Cantelli lemma, applied to a countable sequence $\varepsilon\downarrow0$, we have
\begin{equation}\label{eq:subexp-den}
 \frac{\log Q_{n_j}}{n_j}\rightarrow0
 \quad\text{almost surely}.
\end{equation}
For brevity, write $a_n=-\log\mu_{v,\pi_vX}(C_{\mathbf w_n})$. By \eqref{eq:cond-mass}, \eqref{eq:src-dim} and \eqref{eq:subexp-den}, $a_{n_j}/n_j\rightarrow\dim_H\mu-1$ almost surely. Since the stopped cylinders are nested, $(a_n)$ is nonnegative and nondecreasing. Hence
$$
 \frac{a_{n_j}}{n_{j+1}}\le\frac{a_n}{n}
 \le\frac{a_{n_{j+1}}}{n_j},\qquad n_j\le n\le n_{j+1}.
$$
Since $n_{j+1}/n_j\rightarrow1$, we have
$$
 \frac{-\log\mu_{v,\pi_vX}(C_{\mathbf w_n})}n
 \rightarrow\dim_H\mu-1
 \quad\text{almost surely}.
$$
By \eqref{eq:ball-comp} and disintegration, $\mu_{v,y}$ is exact dimensional of dimension $\dim_H\mu-1$ for $\mu_v$-almost every $y$, hence for $\LL$-almost every $y\in E_v$.
\end{proof}

\subsection{Proof of the set theorem}\label{sec:ssc-app}

\begin{lemma}\label{lem:ssc-app}
In $(\mathrm{I})$, suppose $s_E=\dim_HE>0$.
For every $0<\varepsilon<s_E$, there is an SSC subsystem consisting
of compositions of the maps $S_i$, with attractor
$E_\varepsilon\subset E$ and $\dim_HE_\varepsilon>s_E-\varepsilon$.
Its normalised angles consist of $\theta$ and $\theta+\alpha_i$,
with possible repetitions.

In $(\mathrm{II})$, for every $0<\varepsilon<\dim_HK$, there is an
SSC subsystem consisting of compositions of the maps $S_i$, with
attractor $K_\varepsilon\subset K$ and
$\dim_HK_\varepsilon>\dim_HK-\varepsilon$.
Its contraction logarithms consist of $b_0$ and $b_0+\beta_i$,
with possible repetitions, where $b_0>0$.
The corresponding statement holds for $L$, the maps $S'_j$ and
the logarithms $\beta'_j$.
If the original angle vector satisfies $(D)_\tau$ in $(\mathrm I)$,
or the original joint logarithmic vector satisfies $(D)'_\tau$ in
$(\mathrm{II})$, then the new angle vector in $(\mathrm I)$, or the
new joint logarithmic vector in $(\mathrm{II})$ satisfies the same
condition.
\end{lemma}
\begin{proof}
In $(\mathrm{I})$, remove repeated maps and write
$r_*=\min_i r_i>0$. Fix $s_E-\varepsilon<\sigma<s_E$ and
$\ell_0>2\operatorname{diam}E$. Since $\mathcal H^\sigma(E)=\infty$, a maximal $\ell_0r$-separated
subset of $E$ has at least $c_1r^{-\sigma}$ points for all sufficiently
small $r$, where $c_1>0$ is independent of $r$.
We choose a stopping cylinder at threshold $\log(1/r)$ containing
each point. Each cylinder has diameter at most
$r\operatorname{diam}E$, so the cylinders are pairwise disjoint.
Their contraction ratios belong to $(r_*r,r]$.
Their words have length at most $C_1(1+\log(1/r))$, where $C_1>0$
depends only on the IFS. Hence there are at most
$C_1'(1+\log(1/r))^m$ symbol-count vectors, with $C_1'=(C_1+1)^m$.
We choose a class of maximal size and denote its maps by
$B_{r,1},\ldots,B_{r,M_r}$. Since planar rotations commute, these maps
have a common linear part, and
$$
 M_r\ge\frac{c_1'r^{-\sigma}}{(1+\log(1/r))^m},
$$
where $c_1'=c_1/C_1'$.

Two distinct maps among $\mathrm{Id},S_1,\ldots,S_m$ coincide at
most one point. Since $E$ is infinite, we can choose $x_*\in E$
outside these coincidence points and then a sufficiently small
cylinder $U(E)$ containing $x_*$ such that
\begin{equation}\label{eq:ssc-u}
 U(E),S_1U(E),\ldots,S_mU(E)
 \quad\text{are pairwise disjoint}.
\end{equation}
Fix $U$ independently of $r$, with contraction ratio $r_U$. The system
\begin{equation}\label{eq:ssc-sys}
 \{B_{r,k}\circ U:1\le k\le M_r\}
 \ \cup\
 \{B_{r,1}\circ S_i\circ U:1\le i\le m\}
\end{equation}
satisfies SSC by the separation of the outer cylinders and
\eqref{eq:ssc-u}, and its attractor $E_r$ lies in $E$.
The first $M_r$ maps have a common contraction ratio $\lambda_r$,
with $r_*rr_U<\lambda_r\le rr_U$. Set
$c_1''=c_1'(r_*r_U)^{s_E-\varepsilon}$. Since
$$
 M_r\lambda_r^{s_E-\varepsilon}
 \ge\frac{c_1''r^{s_E-\varepsilon-\sigma}}{(1+\log(1/r))^m}
 \longrightarrow\infty,
$$
the SSC dimension formula \cite[Theorem 9.3]{Fa03} implies
$\dim_HE_r>s_E-\varepsilon$ for all sufficiently small $r$.
The angles in \eqref{eq:ssc-sys} are $\theta_r$ and
$\theta_r+\alpha_i$, where $\theta_r$ is the common angle of
$B_{r,k}\circ U$. Choose one such $r$ and set
$E_\varepsilon=E_r$ and $\theta=\theta_r$.

The same construction applies on the line: the linear parts commute
even when reflections occur, and distinct affine maps coincide at
most one point. Let $b_0$ be the common contraction logarithm of
$B_{r,k}\circ U$. Then $B_{r,1}\circ S_i\circ U$ has contraction
logarithm $b_0+\beta_i$. We apply the construction separately to
$K$ and $L$ to prove $(\mathrm{II})$.

For every integer $k\ne0$, since $k\alpha_i=k(\theta+\alpha_i)-k\theta$, we have
$$
 \max_i\|k\alpha_i\|_{\T}
 \le2\max\left\{\|k\theta\|_{\T},
          \max_i\|k(\theta+\alpha_i)\|_{\T}\right\}.
$$
If the original angle vector satisfies $(D)_\tau$, then so does the
subsystem, with half the original constant. For a product, the new
joint logarithmic vector has entries
$$
 \mathcal B=\{b_0,b'_0\}\cup\{b_0+\beta_i\}_i
                       \cup\{b'_0+\beta'_j\}_j.
$$
For every real $t$, the same subtraction implies
$$
 \max\left\{\max_i\|t\beta_i\|_{\T},
             \max_j\|t\beta'_j\|_{\T}\right\}
 \le2\max_{b\in\mathcal B}\|tb\|_{\T}.
$$
If the original joint logarithmic vector satisfies $(D)'_\tau$, then
the new vector satisfies the same condition, with half the original constant.
\end{proof}

\begin{proof}[Proof of Theorem \ref{thm:sets}]
Write $s_E=\dim_HE>1+2\tau$. By Mattila \cite[Theorem 7.7]{Ma95},
for each fixed $v\in\mathcal D$, we have
\begin{equation}\label{eq:sl-up}
 \dim_H(E\cap\pi_v^{-1}\{y\})\le s_E-1
 \quad\text{for }\LL\text{-almost every }y\in\R.
\end{equation}
Under SSC, the natural self-similar measure, or the product of the natural
self-similar measures in $(\mathrm{II})$, satisfies
$\dim_H\mu=s_E$ and $\mu(B(x,r))\le C_1r^{s_E}$ for $0<r\le1$,
where $C_1>0$ depends only on the IFSs,
see \cite[Theorem 9.3 and its proof]{Fa03}.
Thus $I_s(\mu)<\infty$ for $0<s<s_E$. Choose $1+2\tau<s<s_E$.
By Theorem \ref{thm:ac} and \eqref{eq:sl-up}, we obtain the result
under SSC.

For $0<\varepsilon<s_E-1-2\tau$, Lemma \ref{lem:ssc-app} constructs
an SSC subsystem $E_\varepsilon\subset E$ with
$s_\varepsilon:=\dim_HE_\varepsilon>s_E-\varepsilon$.
In $(\mathrm{II})$, apply the lemma to each factor with errors
whose sum is less than $\varepsilon$, both factors have positive
dimension since their dimensions are at most one and $s_E>1$.
Since Lemma \ref{lem:ssc-app} preserves the original Diophantine condition,
we obtain from the SSC case that
$$
 \inf_{v\in\mathcal V}\LL(E_v)
 \ge\inf_{v\in\mathcal V}\LL(\pi_vE_\varepsilon)>0
$$
for every compact $\mathcal V\subset\mathcal D$.

Fix $v\in\mathcal D$. Since $E_\varepsilon\subset E$, for every
$n\ge1$ and $\mathbf w\in\mathcal P_n$ we have $\Psi_{\mathbf w}(E_\varepsilon)\subset
C_{\mathbf w}=\Psi_{\mathbf w}(E)$.
Choose a compact $\mathcal V_v\subset\mathcal D$ containing all
$v_{n,\mathbf w}$, and let
$$
 \ell_{\varepsilon,v}
 =\inf_{v'\in\mathcal V_v}\LL(\pi_{v'}E_\varepsilon)>0.
$$
For each of the countably many residual vectors $v'$, choose a Borel
set $G_{v'}$ of full Lebesgue measure in $\pi_{v'}E_\varepsilon$
on which the subsystem fibres have dimension $s_\varepsilon-1$.
Set
$$
 G_{\varepsilon,v}
 =\bigcup_{n\ge1}\ \bigcup_{\mathbf w\in\mathcal P_n}
       \bigl(t_{v,\mathbf w}+e^{-n}G_{v_{n,\mathbf w}}\bigr)
 \subset E_v.
$$
By \eqref{eq:st-map}, $\Psi_{\mathbf w}$ maps each fibre of
$E_\varepsilon$ over $G_{v_{n,\mathbf w}}$ into a fibre of $E$.
Since $\Psi_{\mathbf w}$ is bi-Lipschitz
in both cases,
$$
 \dim_H(E\cap\pi_v^{-1}\{y\})\ge s_\varepsilon-1
 \qquad(y\in G_{\varepsilon,v}).
$$
Moreover, for every $\mathbf w\in\mathcal P_n$,
$$
 \LL(G_{\varepsilon,v}\cap\pi_vC_{\mathbf w})
 \ge e^{-n}\ell_{\varepsilon,v}.
$$
By \eqref{eq:cyl-sc} and the argument proving \eqref{eq:den-dev},
$\LL(E_v\setminus G_{\varepsilon,v})=0$.
We let $\varepsilon\downarrow0$ along a countable sequence.
By \eqref{eq:sl-up}, we obtain the required fibre dimension for
$\LL$-almost every $y\in E_v$. The exceptional null set may depend on $v$.
\end{proof}

\appendix

\section{Diophantine criteria and examples}\label{sec:ex}

We now give criteria for checking the Diophantine hypotheses and construct examples satisfying the main theorems. Subsection \ref{sec:dio-prop} relates the integer- and real-frequency conditions and computes the Hausdorff dimension of the exceptional parameter sets. Subsection \ref{sec:arith} provides a rational rank criterion for algebraic vectors. We then use these criteria to construct planar systems with prescribed dimensions in Subsection \ref{sec:pl-ex}, and products and sumsets with prescribed dimensions in Subsection \ref{sec:log-ex}.

\subsection{Integer and real frequencies}\label{sec:dio-prop}

Recall that, for $\bm\alpha=(\alpha_1,\ldots,\alpha_q)\in\R^q$,
$$
 d_{\bm\alpha}(t)=\max_{1\le j\le q}\|t\alpha_j\|_{\T},\qquad t\in\R.
$$

We first relate the Diophantine conditions for integer and real frequencies. Adding the coordinate $1$ allows us to pass between these conditions, while the real-frequency condition is also invariant under nonzero scalar multiplication.

\begin{lemma}\label{lem:real-freq}
Let $\tau>0$.
\begin{enumerate}
\item[(i)] A real vector $\bm\alpha$ satisfies $(D)_\tau$ if and only if $(1,\bm\alpha)$ satisfies $(D)'_\tau$.
\item[(ii)] For every real vector $\bm\beta$ and every real $\lambda\ne0$, $\bm\beta$ satisfies $(D)'_\tau$ if and only if $\lambda\bm\beta$ satisfies $(D)'_\tau$.
\end{enumerate}
\end{lemma}
\begin{proof}
For (i), fix $t\ge1$ and let $k$ be a nearest integer to $t$. Then $1\le k\le3t/2$, and
$$
 d_{\bm\alpha}(k)\le d_{\bm\alpha}(t)+\max_j|\alpha_j|\,\|t\|_{\T}\le(1+\max_j|\alpha_j|)d_{(1,\bm\alpha)}(t).
$$
If $\bm\alpha$ satisfies $(D)_\tau$ with constant $c_1>0$, then \eqref{eq:dio-int} and $k\le3t/2$ imply
$$
 d_{(1,\bm\alpha)}(t)\ge c_1't^{-\tau},\qquad c_1'=\frac{(2/3)^\tau c_1}{1+\max_j|\alpha_j|}.
$$
Conversely, $d_{(1,\bm\alpha)}(k)=d_{\bm\alpha}(k)$ for $k\in\mathbb N$, so \eqref{eq:dio} implies \eqref{eq:dio-int}.

For (ii), suppose that $\bm\beta$ satisfies $(D)'_\tau$ with constant $c_2>0$. By \eqref{eq:dio-small}, for $t\ge1$,
$$
 d_{\lambda\bm\beta}(t)=d_{\bm\beta}(|\lambda|t)\ge 2^{-(\tau+1)}c_2\min\{|\lambda|,|\lambda|^{-\tau}\}t^{-\tau}.
$$
For the reverse implication, we apply the forward implication to $\lambda\bm\beta$ with scalar $1/\lambda$.
\end{proof}

Lemma \ref{lem:real-freq}(ii) applies to $(D)'_\tau$. For angle vectors, we use $(D)_\tau$ in the original normalisation, since this condition is preserved by integer shifts of the angles.

We next show that the integer-frequency condition holds for almost every vector. More precisely, the exceptional vectors have Hausdorff dimension strictly smaller than the ambient dimension.

\begin{proposition}\label{prop:dio-exc}
For every integer $q\ge1$ and every $\tau>1/q$,
$$
 \dim_H\{\bm\alpha\in\R^q:\bm\alpha\text{ fails }(D)_\tau\}=\frac{q+1}{1+\tau}<q.
$$
In particular, $(D)_\tau$ holds for Lebesgue-almost every vector in $\R^q$.
\end{proposition}
\begin{proof}
For $\eta>1/q$, define
$$
 \mathcal J_\eta=\{\bm\alpha\in\R^q:d_{\bm\alpha}(k)<k^{-\eta}\text{ for infinitely many }k\in\mathbb N\}.
$$
By the simultaneous Jarn\'ik-Besicovitch theorem \cite[Theorem 4.12]{BRV16}, with one denominator and $q$ coordinates, we have $\dim_H\mathcal J_\eta=(q+1)/(1+\eta)$ on the unit cube, hence on $\R^q$ by integer translations.

The vectors failing $(D)_\tau$ lie between the sets of vectors admitting approximation at exponents $\tau$ and $\tau+\varepsilon$. To see this, let $\mathcal{E}_\tau=\{\bm\alpha\in\R^q:\inf_{k\ge1}k^\tau d_{\bm\alpha}(k)=0\}$. For every $\varepsilon>0$,
$$
 \mathcal J_{\tau+\varepsilon}\subset\mathcal{E}_\tau\subset\mathcal J_\tau.
$$
For the first inclusion, $k^\tau d_{\bm\alpha}(k)<k^{-\varepsilon}$ along an unbounded sequence. For the second, let $\bm\alpha\in\mathcal{E}_\tau$. If $d_{\bm\alpha}(k_0)=0$ for some $k_0\ge1$, then $d_{\bm\alpha}(jk_0)=0$ for every $j\ge1$. Otherwise the infimum is approached along an unbounded sequence. In either case, $\bm\alpha\in\mathcal J_\tau$. We let $\varepsilon\downarrow0$ to obtain the dimension formula.
\end{proof}

\subsection{Algebraic vectors and rational rank}\label{sec:arith}

For a finite nonzero vector $\bm\beta$, define its reduced rank by
$$
 \delta(\bm\beta)=\dim_{\Q}\operatorname{span}_{\Q}\{\beta_j:j\}-1.
$$
This rank is invariant under nonzero scalar multiplication. We use $(1,\bm\alpha)$ in $(\mathrm I)$ and the contraction logarithms in $(\mathrm{II})$. The hypothesis in Proposition \ref{prop:alg-rank} that a common nonzero scalar multiple has algebraic entries is equivalent to algebraicity of every ratio $\beta_j/\beta_{j_0}$ for one nonzero coordinate $\beta_{j_0}$.

The next proposition shows that, for such algebraic vectors, the reduced rank determines the possible Diophantine exponents. We use Schmidt's theorem for the upper range of exponents and Dirichlet's theorem for the restriction in the other direction.

\begin{proposition}\label{prop:alg-rank}
Let $\bm\beta$ be a finite nonzero real vector whose entries become algebraic after multiplication by a nonzero real scalar, and let $\delta=\delta(\bm\beta)$. If $\delta=0$, then $\bm\beta$ satisfies no $(D)'_\tau$ condition. If $\delta\ge1$, it satisfies $(D)'_\tau$ for every $\tau>1/\delta$, and $(D)'_\tau$ can hold only when $\tau\ge1/\delta$. Consequently, for every $s>1$,
\begin{equation}\label{eq:rank}
 \delta(s-1)>2
 \quad\Longleftrightarrow\quad
 \bm\beta\text{ satisfies }(D)'_\tau
 \text{ for some }0<\tau<(s-1)/2.
\end{equation}
\end{proposition}
\begin{proof}
We reorder the coordinates so that $\beta_1\ne0$ and let $a_j=\beta_j/\beta_1$ and $\mathbf a=(a_j)_j$. Then $\mathbf a$ has algebraic entries, $a_1=1$ and $\delta(\mathbf a)=\delta$. By Lemma \ref{lem:real-freq}(ii), $\mathbf a$ satisfies $(D)'_\tau$ if and only if $\bm\beta$ does. If $\delta=0$, all entries of $\mathbf a$ are rational, so $d_{\mathbf a}$ vanishes at positive integer multiples of a common denominator. Thus $\bm\beta$ satisfies no $(D)'_\tau$ condition.

Suppose that $\delta\ge1$. We first consider a rational basis of the coordinates, since the Diophantine exponent depends on the number of rationally independent coordinates. We choose $1,b_1,\ldots,b_\delta$ among the entries of $\mathbf a$ as a basis for their rational span. Fix $\eta>0$. Since $1,b_1,\ldots,b_\delta$ are linearly independent over $\Q$ and algebraic, Schmidt's theorem on simultaneous approximation \cite[Corollary to Theorem 1]{Sch70} implies that
\begin{equation}\label{eq:alg-app}
 \max_{1\le\ell\le\delta}\|kb_\ell\|_{\T}<k^{-1/\delta-\eta}
\end{equation}
holds for only finitely many $k\in\mathbb N$.

For $\tau>1/\delta$, take $\eta=\tau-1/\delta$. Since $\max_\ell\|kb_\ell\|_{\T}>0$ for every $k\ge1$, we absorb the finitely many exceptions into a constant $c_1>0$, depending only on $\tau$ and $(b_1,\ldots,b_\delta)$, and obtain $(D)_\tau$ for $(b_1,\ldots,b_\delta)$. By Dirichlet's theorem \cite[Theorem 1.2 for $\delta=1$ and Theorem 4.4 with equal weights for $\delta\ge2$]{BRV16}, $\max_\ell\|kb_\ell\|_{\T}\le k^{-1/\delta}$ for arbitrarily large $k$. It follows that no exponent $\tau<1/\delta$ is possible.

We next pass from the basis to the original vector. There are $N\in\mathbb N$ and integers $n_{j\ell}$ such that
$$
 Na_j=n_{j0}+\sum_{\ell=1}^\delta n_{j\ell}b_\ell.
$$
Set $C_2=\max_j\sum_{\ell=0}^\delta|n_{j\ell}|$. By coordinate inclusion and subadditivity of the distance to $\Z$, we have
$$
 d_{(1,b_1,\ldots,b_\delta)}(t)\le d_{\mathbf a}(t),\qquad
 d_{\mathbf a}(Nt)\le C_2d_{(1,b_1,\ldots,b_\delta)}(t),\qquad t\in\R.
$$
If $\mathbf a$ satisfies $(D)'_\tau$ with constant $c_{\mathbf a}>0$, then
$$
 d_{(1,b_1,\ldots,b_\delta)}(t)\ge C_2^{-1}d_{\mathbf a}(Nt)\ge c_{\mathbf a}C_2^{-1}N^{-\tau}t^{-\tau},\qquad t\ge1.
$$
The converse follows from coordinate inclusion. Thus $(D)'_\tau$ for $\mathbf a$ is equivalent to $(D)'_\tau$ for $(1,b_1,\ldots,b_\delta)$, and hence to $(D)_\tau$ for $(b_1,\ldots,b_\delta)$ by Lemma \ref{lem:real-freq}(i). For $\delta\ge1$, the inequalities $1/\delta\le\tau<(s-1)/2$ imply $\delta(s-1)>2$. Conversely, this last condition allows us to choose $1/\delta<\tau<(s-1)/2$. When $\delta=0$, both sides of \eqref{eq:rank} are false.
\end{proof}

\begin{corollary}\label{cor:alg-ac}
In $(\mathrm I)$, suppose that the entries of $\bm\alpha=\bm\alpha(\mathcal S)$ are algebraic, and let $\delta=\delta(1,\bm\alpha)$. In $(\mathrm{II})$, suppose that the ratios between entries of $\bm\beta=\bm\beta(\mathcal S,\mathcal S')$ are algebraic, and let $\delta=\delta(\bm\beta)$. If $I_s(\mu)<\infty$ for some $1<s<2$ with $\delta(s-1)>2$, then the absolute continuity and equivalence conclusions of Theorem \ref{thm:ac} hold, together with its conclusions on the dimensions of conditional measures under the corresponding SSC hypotheses. For every $1<\eta<\delta(s-1)/2$ and compact $\mathcal V\subset\mathcal D$,
$$
 \sup_{v\in\mathcal V}\int f_v(\log(e+f_v))^\eta\dd x<\infty.
$$
If $\delta(\dim_HE-1)>2$, the conclusions of Theorem \ref{thm:sets} hold.
\end{corollary}
\begin{proof}
By Proposition \ref{prop:alg-rank} and Lemma \ref{lem:real-freq}(i), the appropriate Diophantine condition holds for every $\tau>1/\delta$. Fix $1<\eta<\delta(s-1)/2$ and choose
$$
 \frac1\delta<\tau<\frac{s-1}{2\eta}.
$$
Then $(s-1)/(2\tau)>\eta>1$, so we obtain the measure conclusions by Theorem \ref{thm:ac}. For the set conclusions, we choose $1/\delta<\tau<(\dim_HE-1)/2$ and apply Theorem \ref{thm:sets}.
\end{proof}

We next give algebraic vectors for which the endpoint exponent $\tau=1/d$ also holds. This will allow us to use the endpoint in the examples below.

For an integer $d\ge1$, let $a_*=2^{1/(d+1)}$. By Eisenstein's criterion with $p=2$, the polynomial $x^{d+1}-2$ is irreducible over $\Q$. Thus $a_*$ has degree $d+1$ over $\Q$, so $1,a_*,\ldots,a_*^d$ are linearly independent over $\Q$ and
\begin{equation}\label{eq:pow-rank}
 \delta(1,a_*,\ldots,a_*^d)=d.
\end{equation}

The following lemma is a special case of Perron's theorem on algebraic badly approximable vectors \cite{Per21}. We include an elementary proof for completeness. For further explicit algebraic examples see e.g. \cite{Cus83,Bri03,Bug15}.

\begin{lemma}\label{lem:pow-end}
Let $d\ge1$ be an integer, $a_*=2^{1/(d+1)}$ and $\mathbf a=(1,a_*,\ldots,a_*^d)$. Then $(a_*,\ldots,a_*^d)$ satisfies $(D)_{1/d}$, and $\mathbf a$ satisfies $(D)'_{1/d}$.
\end{lemma}
\begin{proof}
We use an integer matrix having $\mathbf a$ as an eigenvector. The determinant of the iterates of an integer vector will give a lower bound for its distance from the line spanned by $\mathbf a$.

Let $T$ be the integer matrix defined by
$$
 T(x_0,\ldots,x_d)=(x_1,\ldots,x_d,2x_0).
$$
Then $T\mathbf a=a_*\mathbf a$ and $T^{d+1}=2I$. For every nonzero integer vector $\mathbf n$, we have
$$
 \det(\mathbf n,T\mathbf n,\ldots,T^d\mathbf n)\ne0.
$$
Otherwise the integer columns are rationally dependent, so $P(T)\mathbf n=0$ for a nonzero $P\in\Q[x]$ of degree at most $d$. Since $x^{d+1}-2$ is irreducible and $\deg P\le d$, these polynomials are coprime. By B\'ezout's identity and $T^{d+1}=2I$, the matrix $P(T)$ is invertible, contradicting $P(T)\mathbf n=0$.

We now choose $\mathbf n$ close to $t\mathbf a$. The determinant is a nonzero integer, while the column operations below bound it in terms of the approximation error.

Fix $t\ge1$, choose an integer vector $\mathbf n$ nearest to $t\mathbf a$ coordinatewise, and write
$$
 \mathbf e=\mathbf n-t\mathbf a
$$
so that $d_{\mathbf a}(t)=\|\mathbf e\|_\infty$. Since $t\ge1$, the first coordinate of $\mathbf n$ is nonzero. There is $C_d>0$, depending only on $d$, such that $\|\mathbf n\|_\infty\le C_dt$. Also $\|\mathbf e\|_\infty\le1/2$. For each $1\le j\le d$, replace the column $T^j\mathbf n$ by $T^j\mathbf n-a_*^j\mathbf n$. These column operations do not change the determinant. Since
$$
 T^j\mathbf n-a_*^j\mathbf n=(T^j-a_*^jI)\mathbf e,
$$
Hadamard's inequality and the nonzero integer determinant imply
$$
 1\le|\det(\mathbf n,T\mathbf n,\ldots,T^d\mathbf n)|
 \le\|\mathbf n\|_2\prod_{j=1}^d\|(T^j-a_*^jI)\mathbf e\|_2
 \le C_d't\,d_{\mathbf a}(t)^d.
$$
Here $C_d'>0$ depends only on $d$. Thus, with $c_d=(C_d')^{-1/d}$, we have $d_{\mathbf a}(t)\ge c_d t^{-1/d}$ for every real $t\ge1$, which proves $(D)'_{1/d}$. At integer $t$, the initial coordinate $1$ contributes zero, so $(a_*,\ldots,a_*^d)$ satisfies $(D)_{1/d}$.
\end{proof}

By Lemmas \ref{lem:pow-end} and \ref{lem:real-freq}(ii), we may take $\tau=1/d$ for $(a_*,\ldots,a_*^d)$ as angles and for positive scalar multiples of $\mathbf a$ as contraction logarithms. If $I_s(\mu)<\infty$ for some $1<s<2$ with $d(s-1)\ge2$, Theorem \ref{thm:ac} applies with logarithmic moment exponent $d(s-1)/2$.

\subsection{Planar systems of prescribed dimension}\label{sec:pl-ex}

We now construct planar self-similar systems with any prescribed dimension $s_E\in(1,2)$. We use a grid of translations to ensure SSC and algebraic rotation angles to obtain the required Diophantine condition.

Fix $s_E\in(1,2)$ and choose an integer $q\ge3$ such that $q^2(s_E-1)>2$ and
\begin{equation}\label{eq:sep}
 0<\rho:=q^{-2/s_E}<\frac1{1+\sqrt2(q-1)}.
\end{equation}
Such $q$ exists since $2/s_E>1$. Let $d=q^2$ and $a_j=2^{j/(d+1)}$, $1\le j\le d$. By \eqref{eq:pow-rank}, $1,a_1,\ldots,a_d$ are linearly independent over $\Q$. For $0\le i,j<q$, define
$$
 t_{ij}=\frac{(2i-q+1,\,2j-q+1)}{2(q-1)},\qquad
 S_{ij}(x)=\rho O_{a_{qi+j+1}}x+t_{ij}.
$$
Since $|t_{ij}|\le1/\sqrt2$, the maps preserve the closed disk of radius $r_\rho=1/(\sqrt2(1-\rho))$. Their image disks have radii $\rho r_\rho$ and centre separation at least $1/(q-1)$. By \eqref{eq:sep}, $2\rho r_\rho<1/(q-1)$, so the attractor $E$ satisfies SSC. Since $q^2\rho^{s_E}=1$, by the SSC dimension formula we have $\dim_HE=s_E$.

We now apply the algebraic Diophantine criterion to the rotation angles. The self-similar measure $\mu$ with weights $q^{-2}=\rho^{s_E}$ is $s_E$-Ahlfors regular, so $I_s(\mu)<\infty$ for $0<s<s_E$. For $1<\eta<d(s_E-1)/2$, choose
$
 1+\frac{2\eta}{d}<s<s_E.
$
By Corollary \ref{cor:alg-ac}, every nonzero linear image of $\mu$ is equivalent to Lebesgue measure on its support and satisfies the logarithmic moment bound with exponent $\eta$, uniformly for $v$ in compact subsets of $\R^2\setminus\{0\}$. Every projection of $E$ has positive length, and for each fixed $v\ne0$,
$$
 \dim_H\bigl(E\cap\pi_v^{-1}\{y\}\bigr)=s_E-1
 \quad\text{for }\LL\text{-almost every }y\in\pi_vE.
$$

\subsection{Products and sumsets of prescribed dimensions}\label{sec:log-ex}

We first construct one-dimensional SSC systems with prescribed dimension and contraction logarithms proportional to a chosen algebraic vector. We then combine two such systems so that their joint logarithmic vector satisfies the required Diophantine condition.

Let $\mathbf a=(a_1,\ldots,a_m)$ have positive algebraic entries, with $m\ge2$, and fix $s_K\in(0,1)$. Since $\sum_i e^{-s_Kca_i}$ decreases continuously from $m$ to zero as $c$ increases from zero to infinity, there is a unique $c>0$ satisfying $\sum_i e^{-s_Kca_i}=1$. Define
\begin{equation}\label{eq:log-fam}
 r_i=e^{-ca_i},\qquad 1\le i\le m.
\end{equation}
Since $s_K<1$, we have $\sum_i r_i<\sum_i r_i^{s_K}=1$. Let $S_i(x)=r_ix+d_i$, where
\begin{equation}\label{eq:trans}
 d_i=\sum_{h<i}r_h+
 \frac{i-1}{m-1}\left(1-\sum_{h=1}^m r_h\right),
 \qquad 1\le i\le m.
\end{equation}
The intervals $S_i([0,1])$ have equal positive gaps $(1-\sum_h r_h)/(m-1)$, with leftmost endpoint $0$ and rightmost endpoint $1$. Hence the attractor $K$ satisfies SSC and $\dim_HK=s_K$. If $\delta=\delta(\mathbf a)\ge1$, its logarithmic vector $c\mathbf a$ satisfies $(D)'_\tau$ for every $\tau>1/\delta$ by Proposition \ref{prop:alg-rank}. The exponent range is independent of $c$, while the Diophantine constant may depend on $c$.

We now construct examples using logarithms from both IFSs. The main point is that the Diophantine condition can be obtained from the first logarithmic vector together with differences of logarithms from the second. This allows us to choose the dimensions of the two sets while keeping the required rational rank.

Fix $m,n\ge2$ and $s_K,s_L\in(0,1)$ with
$
 \delta(s_K+s_L-1)>2$ with $\delta=m+n-2.
$
Choose $\mathbf a=(a_1,\ldots,a_{m+n-1})$ with positive algebraic entries linearly independent over $\Q$, for example,
$
 a_\ell=2^{(\ell-1)/(m+n-1)}$, $1\le\ell\le m+n-1.
$
Choose $c>0$ so that $\sum_{i=1}^m e^{-s_Kca_i}=1$, and define
$$
 b_0=\frac1{s_L}\log\left(1+
 \sum_{j=2}^n e^{-s_Lca_{m+j-1}}\right),\qquad
 r_i=e^{-ca_i},
$$
$$
 r'_1=e^{-b_0},\qquad
 r'_j=e^{-(b_0+ca_{m+j-1})},\quad 2\le j\le n.
$$
Then
$$
 \sum_{j=1}^n (r'_j)^{s_L}
 =e^{-s_Lb_0}\left(1+\sum_{j=2}^n e^{-s_Lca_{m+j-1}}\right)=1.
$$
Since $s_K,s_L<1$, we have $\sum_i r_i<1$ and $\sum_j r'_j<1$. By \eqref{eq:trans}, these ratios define SSC systems with attractors $K,L\subset[0,1]$ and dimensions $s_K,s_L$.

The common term $b_0$ disappears when taking differences of the logarithms of the second system. We can therefore verify the Diophantine condition using the algebraic vector $\mathbf a$, without requiring $b_0$ itself to be algebraic.

Write $\beta_i=\log(1/r_i)$ and $\beta'_j=\log(1/r'_j)$. The vector
$$
 \bm\gamma=(\beta_1,\ldots,\beta_m,\,
 \beta'_2-\beta'_1,\ldots,\beta'_n-\beta'_1)=c\mathbf a
$$
has reduced rank $\delta=m+n-2$, so Proposition \ref{prop:alg-rank} implies $(D)'_\tau$ for every $\tau>1/\delta$. For $a_\ell=2^{(\ell-1)/(m+n-1)}$, Lemmas \ref{lem:pow-end} and \ref{lem:real-freq}(ii) imply $(D)'_{1/\delta}$ for $\bm\gamma=c\mathbf a$. Since
$$
 d_{\bm\gamma}(t)\le2d_{\bm\beta(\mathcal S,\mathcal S')}(t),
 \qquad t\in\R,
$$
the full logarithmic vector satisfies \eqref{eq:sum-dio} for every $\tau>1/\delta$, and also for $\tau=1/\delta$ in this power family. Neither Diophantine conclusion requires $b_0/c$ to be algebraic.

We can now apply the main theorems. Choose $1/\delta<\tau<(s_K+s_L-1)/2$. By Theorem \ref{thm:sums}, $\LL(uK+L)>0$ for all $u\ne0$, uniformly on compact subsets of $\R\setminus\{0\}$, and, for each fixed $u\ne0$,
$$
 \dim_H\bigl(L\cap(y-uK)\bigr)=s_K+s_L-1
 \quad\text{for }\LL\text{-almost every }y\in uK+L.
$$
The natural measures $\nu,\nu'$ with weights $r_i^{s_K}$ and $(r'_j)^{s_L}$ are Ahlfors regular. Their product has finite $s$-energy for $0<s<s_K+s_L$. We choose $1+2\tau<s<s_K+s_L$ and apply Theorem \ref{thm:ac}. Thus $\nu_u*\nu'\sim\LL|_{uK+L}$ for every $u\ne0$, with the logarithmic moment bound \eqref{eq:log-mom} uniform on compact subsets of $\R\setminus\{0\}$.

The joint condition can hold even when neither factor alone has a sufficiently strong Diophantine exponent. For $m=2$, $n=3$ and $s_K+s_L>5/3$, the difference vector has reduced rank $\delta=3$, so every $1/3<\tau<(s_K+s_L-1)/2$ is admissible. For each factor, we divide its logarithmic vector by one coordinate. By Lemma \ref{lem:real-freq}(i)-(ii) and Dirichlet's theorem in one or two dimensions, a vector with two logarithms can satisfy $(D)'_\tau$ only if $\tau\ge1$, and one with three logarithms only if $\tau\ge1/2$. Since $(s_K+s_L-1)/2<1/2$, neither factor alone suffices for this application of Theorem \ref{thm:sums}.

Finally, we construct a fixed SSC set of dimension arbitrarily close to $1$ with a correspondingly small Diophantine exponent. For $\varepsilon\in(0,1)$, let $s_K=1-\varepsilon/2$. Choose an integer $m>1+4/\varepsilon$ and let
$$
 a_i=2^{(i-1)/m},\qquad 1\le i\le m.
$$
By Lemma \ref{lem:pow-end}, this vector satisfies $(D)'_{\varepsilon/4}$ since $1/(m-1)<\varepsilon/4$. By \eqref{eq:log-fam}-\eqref{eq:trans}, we have a fixed SSC set $K$ of dimension $1-\varepsilon/2$ whose logarithmic vector satisfies $(D)'_{\varepsilon/4}$.

\section*{Acknowledgements}

We thank Amir Algom, Simon Baker, Pablo Shmerkin and P\'eter Varj\'u for helpful discussions.

\section*{AI statement}

AI tools were used to assist with language editing, cross-referencing and proof-reading in the final stages of the project. The mathematical ideas, arguments and proofs are the authors' own. The authors take responsibility for the content of the manuscript.


\begin{thebibliography}{99}
\bibitem{Alg20}
A. Algom,
\emph{Slicing theorems and rigidity phenomena for self-affine carpets},
Proc. Lond. Math. Soc. (3) \textbf{121} (2020), no. 2, 312-353.
\href{https://doi.org/10.1112/plms.12325}{doi:10.1112/plms.12325}.

\bibitem{ARW26b}
A. Algom, F. Rodriguez Hertz and Z. Wang,
\emph{Absolutely Continuous Convolutions and Projections of Fractal Measures},
preprint, 2026.
\href{https://arxiv.org/abs/2609.36035}{arXiv:2609.36035}.

\bibitem{AW23}
A. Algom and M. Wu,
\emph{Improved versions of some Furstenberg type slicing theorems for self-affine carpets},
Int. Math. Res. Not. IMRN \textbf{2023} (2023), no. 3, 2304-2343.
\href{https://doi.org/10.1093/imrn/rnab318}{doi:10.1093/imrn/rnab318}.

\bibitem{BKK26}
B. B\'ar\'any, A. K\"aenm\"aki and I. Kolossv\'ary,
\emph{Projections of self-affine sets onto lines}, preprint, 2026.
\href{https://arxiv.org/abs/2607.14740}{arXiv:2607.14740}.

\bibitem{BKY26}
B. B\'ar\'any, A. K\"aenm\"aki and H. Yu,
\emph{Finer geometry of planar self-affine sets},
Adv. Math. \textbf{493} (2026), Paper No. 110914.
\href{https://doi.org/10.1016/j.aim.2026.110914}{doi:10.1016/j.aim.2026.110914}.

\bibitem{BS16}
Y. Benoist and N. de Saxc\'e,
\emph{A spectral gap theorem in simple Lie groups},
Invent. Math. \textbf{205} (2016), no. 2, 337-361.
\href{https://doi.org/10.1007/s00222-015-0636-2}{doi:10.1007/s00222-015-0636-2}.

\bibitem{BRV16}
V. Beresnevich, F. Ram\'irez and S. Velani,
\emph{Metric Diophantine approximation: aspects of recent work},
in \emph{Dynamics and Analytic Number Theory},
D. Badziahin, A. Gorodnik and N. Peyerimhoff (eds.),
London Mathematical Society Lecture Note Series, vol. 437,
Cambridge University Press, Cambridge, 2016, pp. 1-95.
\href{https://doi.org/10.1017/9781316402696.002}{doi:10.1017/9781316402696.002}.

\bibitem{Bor21}
B. Borda, \emph{Equidistribution of random walks on compact
groups II. The Wasserstein metric},
Bernoulli \textbf{27} (2021), 2598-2623.
\href{https://doi.org/10.3150/21-BEJ1324}{doi:10.3150/21-BEJ1324}.

\bibitem{Bou10}
J. Bourgain, \emph{The discretized sum-product and projection theorems},
J. Anal. Math. \textbf{112} (2010), 193-236.
\href{https://doi.org/10.1007/s11854-010-0028-x}{doi:10.1007/s11854-010-0028-x}.

\bibitem{BG08}
J. Bourgain and A. Gamburd,
\emph{On the spectral gap for finitely-generated subgroups of $\mathrm{SU}(2)$},
Invent. Math. \textbf{171} (2008), no. 1, 83-121.
\href{https://doi.org/10.1007/s00222-007-0072-z}{doi:10.1007/s00222-007-0072-z}.

\bibitem{Bri03}
K. Briggs, \emph{Some explicit badly approximable pairs},
J. Number Theory \textbf{103} (2003), no. 1, 71-76.
\href{https://doi.org/10.1016/S0022-314X(03)00104-5}{doi:10.1016/S0022-314X(03)00104-5}.

\bibitem{BJ19}
C. Bruce and X. Jin,
\emph{Projections of Gibbs measures on self-conformal sets},
Nonlinearity \textbf{32} (2019), no. 2, 603-621.
\href{https://doi.org/10.1088/1361-6544/aaec9f}{doi:10.1088/1361-6544/aaec9f}.
Corrected version addressing an error in the definition of faithful maps:
\href{https://arxiv.org/abs/1801.06468v3}{arXiv:1801.06468v3}.

\bibitem{Bug15}
Y. Bugeaud, \emph{On the multiples of a badly approximable vector},
Acta Arith. \textbf{168} (2015), no. 1, 71-81.
\href{https://doi.org/10.4064/aa168-1-4}{doi:10.4064/aa168-1-4}.

\bibitem{Bur73}
D. L. Burkholder,
\emph{Distribution Function Inequalities for Martingales},
Ann. Probab. \textbf{1} (1973), no. 1, 19-42.
\href{https://doi.org/10.1214/aop/1176997023}{doi:10.1214/aop/1176997023}.

\bibitem{BG70}
D. L. Burkholder and R. F. Gundy,
\emph{Extrapolation and interpolation of quasi-linear operators on martingales},
Acta Math. \textbf{124} (1970), 249-304.
\href{https://doi.org/10.1007/BF02394573}{doi:10.1007/BF02394573}.

\bibitem{Cus83}
T. W. Cusick,
\emph{The two-dimensional Diophantine approximation constant. II},
Pacific J. Math. \textbf{105} (1983), no. 1, 53-67.
\href{https://doi.org/10.2140/pjm.1983.105.53}{doi:10.2140/pjm.1983.105.53}.

\bibitem{Ess66}
C. G. Esseen, \emph{On the Kolmogorov-Rogozin inequality for the
concentration function},
Z. Wahrscheinlichkeitstheorie verw. Gebiete \textbf{5} (1966), 210-216.
\href{https://doi.org/10.1007/BF00533057}{doi:10.1007/BF00533057}.

\bibitem{Ero07}
K. I. Ero\u{g}lu,
\emph{On planar self-similar sets with a dense set of rotations},
Ann. Acad. Sci. Fenn. Math. \textbf{32} (2007), no. 2, 409-424.
\href{https://arxiv.org/abs/math/0603181}{arXiv:math/0603181}.

\bibitem{Fa03}
K. Falconer, \emph{Fractal Geometry: Mathematical Foundations
and Applications}, second edition, John Wiley \& Sons, 2003.
\href{https://doi.org/10.1002/0470013850}{doi:10.1002/0470013850}.

\bibitem{Fa82}
K. J. Falconer, \emph{Hausdorff dimension and the exceptional set of projections},
Mathematika \textbf{29} (1982), no. 1, 109-115.
\href{https://doi.org/10.1112/S0025579300012201}{doi:10.1112/S0025579300012201}.

\bibitem{Fa26}
K. J. Falconer, \emph{Seventy Years of Fractal Projections},
preprint, 2026.
\href{https://arxiv.org/abs/2602.22002}{arXiv:2602.22002}.

\bibitem{FFJ15}
K. J. Falconer, J. M. Fraser and X. Jin,
\emph{Sixty years of fractal projections},
in \emph{Fractal Geometry and Stochastics V},
Progr. Probab. \textbf{70}, Birkh\"auser, 2015, 3-25.
\href{https://doi.org/10.1007/978-3-319-18660-3_1}{doi:10.1007/978-3-319-18660-3\_1}.

\bibitem{FJ14}
K. Falconer and X. Jin,
\emph{Exact dimensionality and projections of random self-similar
measures and sets},
J. Lond. Math. Soc. (2) \textbf{90} (2014), no. 2, 388-412.
\href{https://doi.org/10.1112/jlms/jdu031}{doi:10.1112/jlms/jdu031}.

\bibitem{FJ15}
K. Falconer and X. Jin, \emph{Dimension conservation for self-similar
sets and fractal percolation},
Int. Math. Res. Not. IMRN \textbf{2015} (2015), no. 24, 13260-13289.
\href{https://doi.org/10.1093/imrn/rnv103}{doi:10.1093/imrn/rnv103}.

\bibitem{FJS10}
A. Ferguson, T. Jordan and P. Shmerkin,
\emph{The Hausdorff dimension of the projections of self-affine carpets},
Fund. Math. \textbf{209} (2010), no. 3, 193-213.
\href{https://doi.org/10.4064/fm209-3-1}{doi:10.4064/fm209-3-1}.

\bibitem{FS18}
J. M. Fraser and T. Sahlsten,
\emph{On the Fourier analytic structure of the Brownian graph},
Anal. PDE \textbf{11} (2018), no. 1, 115--132.
\href{https://doi.org/10.2140/apde.2018.11.115}{doi:10.2140/apde.2018.11.115}.

\bibitem{LG16}
J.-F. Le Gall, \emph{Brownian Motion, Martingales, and Stochastic Calculus}, Graduate Texts in Mathematics, vol. 274, Springer, 2016.
\href{https://doi.org/10.1007/978-3-319-31089-3}{doi:10.1007/978-3-319-31089-3}.

\bibitem{HS12}
M. Hochman and P. Shmerkin,
\emph{Local entropy averages and projections of fractal measures},
Ann. of Math. (2) \textbf{175} (2012), no. 3, 1001-1059.
\href{https://doi.org/10.4007/annals.2012.175.3.1}{doi:10.4007/annals.2012.175.3.1}.

\bibitem{HKY93}
B. R. Hunt, I. Kan and J. A. Yorke,
\emph{When Cantor sets intersect thickly},
Trans. Amer. Math. Soc. \textbf{339} (1993), no. 2, 869-888.
\href{https://doi.org/10.1090/S0002-9947-1993-1117219-8}{doi:10.1090/S0002-9947-1993-1117219-8}.

\bibitem{JS26}
X. Jin and T. Sahlsten,
\emph{Absolute continuity of projections and dimension conservation
for self-similar sets and measures}, preprint, 2026.
\href{https://arxiv.org/abs/2609.29301v2}{arXiv:2609.29301v2}.

\bibitem{KO23}
A. K\"aenm\"aki and T. Orponen,
\emph{Absolute continuity in families of parametrised non-homogeneous
self-similar measures},
J. Fractal Geom. \textbf{10} (2023), no. 1/2, 169-207.
\href{https://doi.org/10.4171/JFG/127}{doi:10.4171/JFG/127}.

\bibitem{Ka68}
R. Kaufman, \emph{On Hausdorff dimension of projections},
Mathematika \textbf{15} (1968), no. 2, 153-155.
\href{https://doi.org/10.1112/S0025579300002503}{doi:10.1112/S0025579300002503}.

\bibitem{LS22}
J. Li and T. Sahlsten, \emph{Trigonometric series and self-similar sets},
J. Eur. Math. Soc. \textbf{24} (2022), no. 1, 341-368.
\href{https://doi.org/10.4171/JEMS/1102}{doi:10.4171/JEMS/1102}.

\bibitem{Mar54}
J. M. Marstrand, \emph{Some fundamental geometrical properties of plane sets
of fractional dimensions},
Proc. London Math. Soc. (3) \textbf{4} (1954), no. 1, 257-302.
\href{https://doi.org/10.1112/plms/s3-4.1.257}{doi:10.1112/plms/s3-4.1.257}.

\bibitem{Ma75}
P. Mattila, \emph{Hausdorff dimension, orthogonal projections and intersections with planes},
Ann. Acad. Sci. Fenn. Ser. A I Math. \textbf{1} (1975), no. 2, 227-244.
\href{https://doi.org/10.5186/aasfm.1975.0110}{doi:10.5186/aasfm.1975.0110}.

\bibitem{Ma95}
P. Mattila, \emph{Geometry of Sets and Measures in Euclidean Spaces:
Fractals and Rectifiability},
Cambridge Studies in Advanced Mathematics, vol. 44,
Cambridge University Press, Cambridge, 1995.
\href{https://doi.org/10.1017/CBO9780511623813}{doi:10.1017/CBO9780511623813}.

\bibitem{Ob12}
D. M. Oberlin, \emph{Restricted Radon transforms and projections of planar sets},
Canad. Math. Bull. \textbf{55} (2012), no. 4, 815-820.
\href{https://doi.org/10.4153/CMB-2011-064-6}{doi:10.4153/CMB-2011-064-6}.

\bibitem{OS23}
T. Orponen and P. Shmerkin,
\emph{On the Hausdorff dimension of Furstenberg sets and orthogonal
projections in the plane},
Duke Math. J. \textbf{172} (2023), no. 18, 3559-3632.
\href{https://doi.org/10.1215/00127094-2022-0103}{doi:10.1215/00127094-2022-0103}.

\bibitem{PS00}
Y. Peres and W. Schlag,
\emph{Smoothness of projections, Bernoulli convolutions, and the dimension
of exceptions},
Duke Math. J. \textbf{102} (2000), no. 2, 193-251.
\href{https://doi.org/10.1215/S0012-7094-00-10222-0}{doi:10.1215/S0012-7094-00-10222-0}.

\bibitem{PS09}
Y. Peres and P. Shmerkin, \emph{Resonance between Cantor sets},
Ergodic Theory Dynam. Systems \textbf{29} (2009), 201-221.
\href{https://arxiv.org/abs/0705.2628}{arXiv:0705.2628}.

\bibitem{Per21}
O. Perron, \emph{\"Uber diophantische Approximationen},
Math. Ann. \textbf{83} (1921), 77-84.
\href{https://doi.org/10.1007/BF01464229}{doi:10.1007/BF01464229}.

\bibitem{Py25}
A. Py\"or\"al\"a,
\emph{On dissonance and orthogonal projections of self-conformal measures},
preprint, 2025, revised 2026.
\href{https://arxiv.org/abs/2511.14493v2}{arXiv:2511.14493v2}.

\bibitem{Py25b}
A. Py\"or\"al\"a,
\emph{The dimension of projections of planar diagonal self-affine measures},
Ann. Fenn. Math. \textbf{50} (2025), no. 1, 59-78.
\href{https://doi.org/10.54330/afm.156815}{doi:10.54330/afm.156815}.

\bibitem{Ra17}
A. Rapaport,
\emph{A self-similar measure with dense rotations, singular
projections and discrete slices},
Adv. Math. \textbf{321} (2017), 529-546.
\href{https://arxiv.org/abs/1703.09977}{arXiv:1703.09977}.

\bibitem{Ra20}
A. Rapaport,
\emph{On self-similar measures with absolutely continuous projections
and dimension conservation in each direction},
Ergodic Theory Dynam. Systems \textbf{40} (2020), no. 12, 3438-3456.
\href{https://doi.org/10.1017/etds.2019.39}{doi:10.1017/etds.2019.39}.

\bibitem{RW23}
K. Ren and H. Wang, \emph{Furstenberg sets estimate in the plane},
preprint, 2023.
\href{https://arxiv.org/abs/2308.08819v3}{arXiv:2308.08819v3}.

\bibitem{Sch70}
W. M. Schmidt,
\emph{Simultaneous approximation to algebraic numbers by rationals},
Acta Math. \textbf{125} (1970), 189-201.
\href{https://doi.org/10.1007/BF02392334}{doi:10.1007/BF02392334}.

\bibitem{Sh14}
P. Shmerkin, \emph{On the exceptional set for absolute continuity
of Bernoulli convolutions},
Geom. Funct. Anal. \textbf{24} (2014), 946-958.
\href{https://arxiv.org/abs/1303.3992}{arXiv:1303.3992}.

\bibitem{Sh19a}
P. Shmerkin, \emph{On Furstenberg's intersection conjecture,
self-similar measures, and the $L^q$ norms of convolutions},
Ann. of Math. (2) \textbf{189} (2019), no. 2, 319-391.
\href{https://doi.org/10.4007/annals.2019.189.2.1}{doi:10.4007/annals.2019.189.2.1}.

\bibitem{SS16}
P. Shmerkin and B. Solomyak, \emph{Absolute continuity of self-similar
measures, their projections and convolutions},
Trans. Amer. Math. Soc. \textbf{368} (2016), 5125-5151.
\href{https://arxiv.org/abs/1406.0204}{arXiv:1406.0204}.

\bibitem{SU26}
P. Shmerkin and S. Usuki,
\emph{Markov-driven self-similar measures and non-homogeneous slicing},
preprint, 2026.
\href{https://arxiv.org/abs/2611.00962}{arXiv:2611.00962}.

\bibitem{Va13}
P. P. Varj\'u, \emph{Random walks in compact groups},
Doc. Math. \textbf{18} (2013), 1137-1175.
\href{https://doi.org/10.4171/DM/423}{doi:10.4171/DM/423}.

\bibitem{Wu19}
M. Wu, \emph{A proof of Furstenberg's conjecture on the intersections
of $\times p$- and $\times q$-invariant sets},
Ann. of Math. (2) \textbf{189} (2019), no. 3, 707-751.
\href{https://doi.org/10.4007/annals.2019.189.3.2}{doi:10.4007/annals.2019.189.3.2}.

\end{thebibliography}
\end{document}